\pdfoutput=1

\documentclass{article}
\usepackage{a4wide}
\usepackage[english]{babel}
\usepackage{amsmath}
\usepackage{amssymb}
\usepackage{siunitx}
\usepackage{amsfonts}
\usepackage{leftidx}
\usepackage{graphicx}
\usepackage{subcaption}
\usepackage[lined,boxed,commentsnumbered]{algorithm2e}

\newcommand{\citep}{\cite}

\numberwithin{equation}{section}
\numberwithin{figure}{section}

\begin{document}

\title{Condition number analysis and preconditioning of the finite cell method}
\author{F. de Prenter\footnote{Corresponding author, tel.: +31 61 516 2599, email: f.d.prenter@tue.nl}, C.V. Verhoosel, G.J. van Zwieten, E.H. van Brummelen \\ Eindhoven University of Technology}
\date{\today}

\maketitle

\begin{abstract}
\noindent The (Isogeometric) Finite Cell Method -- in which a domain is immersed in a structured background mesh -- suffers from conditioning problems when cells with small volume fractions occur. In this contribution, we establish a rigorous scaling relation between the condition number of (I)FCM system matrices and the smallest cell volume fraction. Ill-conditioning stems either from basis functions being small on cells with small volume fractions, or from basis functions being nearly linearly dependent on such cells. Based on these two sources of ill-conditioning, an algebraic preconditioning technique is developed, which is referred to as \emph{Symmetric Incomplete Permuted Inverse Cholesky} (SIPIC). A detailed numerical investigation of the effectivity of the SIPIC preconditioner in improving (I)FCM condition numbers and in improving the convergence speed and accuracy of iterative solvers is presented for the Poisson problem and for two- and three-dimensional problems in linear elasticity, in which Nitche's method is applied in either the normal or tangential direction. The accuracy of the preconditioned iterative solver enables mesh convergence studies of the finite cell method.
\end{abstract}


\section{Introduction}
The Isogeometric Finite Cell Method (IFCM, \emph{e.g.,} \cite{Schillinger2012,Rank2012,Schillinger2014})
combines Isogeometric Analysis (IGA, \citep{Hughes2005,Cottrell2009}) with the Finite Cell Method (FCM, \cite{Parvizian2007,Duester2008}).
This combination enables application of IGA to trimmed, coupled and overlapping domains in Computer Aided Design (CAD) \cite{Schmidt2012,Ruess2013,Ruess2014},
topologically complex structures such as porous materials, composites and scanned data \cite{Ruess2012,Verhoosel2015},
and problems with moving boundaries such as Fluid Structure Interactions \cite{Burman2014,Kamensky2015,Hsu2015}, without laborious (re-)meshing procedures.

While FCM was initially introduced as a $p$-FEM method \cite{Parvizian2007}, it is nowadays commonly used in combination with IGA.
The main advantage of FCM over standard finite element methods is that the mesh does not need to match the boundaries of the physical domain.
Instead, the physical domain is immersed in a topologically and geometrically simpler encapsulating mesh, on which a myriad of approximation spaces -- including Non-Uniform Rational B-splines (NURBS) -- can be formed. In (I)FCM, the complexity of the physical domain is captured by an advanced integration procedure for the cells that intersect the domain boundary. Dirichlet boundary conditions on the immersed boundaries are weakly imposed using Nitsche's method \citep{Nitsche}.

When the encapsulating mesh contains cells that only intersect the physical domain on a small fraction of their volume, FCM has been found to be prone to conditioning problems. These conditioning problems impede solving the resulting system of equations. To phrase the recent review article by Schillinger and Ruess \citep{Schillinger2014}:
``$\ldots$ the development of suitable preconditioning techniques that open the door for efficient iterative solution methods in the finite cell method seems very desirable $\ldots$''.
As a matter of fact, various conditioning strategies have already been proposed. The most prominent of these are:
\begin{itemize}
\item \emph{Fictitious domain stiffness:} In order to increase the contribution of basis functions of which only a small fraction of their support intersects the physical domain, many authors have used a virtual stiffness, \emph{e.g.,} \cite{Schillinger2014}. This implies that volumetric operators are not only integrated over the physical domain but also -- multiplied by a small parameter -- over the fictitious domain.
See \cite{RankMath} for a mathematical analysis of this approach.
\item \emph{Basis function manipulation:} The adverse effect of basis functions with small supports within the physical domain can be ameliorated by eliminating them from the system. This can either be done by combining them with geometrically nearby functions \cite{Rueberg2012,Rueberg2014} (as is also done in Web-splines \cite{Hoellig2002,Hoellig2005}) or by simply excluding them from the approximation space, \emph{e.g.,} \citep{Embar2010,Sanches2011,Verhoosel2015}. 
\item \emph{Ghost penalty:} The contribution of basis function with small supports within the physical domain can also be increased by adding an extra term to the variational formulation. An example of such a modification of the formulation is the addition of the Ghost penalty term \citep{GhostPenalty}.
\item \emph{Basis function scaling:} In the context of linear basis functions and the extended finite element method (XFEM),
diagonal scaling of the system matrix has been demonstrated to resolve conditioning problems in specific cases \cite{DiagonalXFEM}.
\end{itemize}
Although these methods (and combinations thereof) can be effective in improving the conditioning of finite cell systems,
they generally manipulate the weak formulation and/or approximation space, which inevitably affects the solution and may compromise stability properties of the employed approximation spaces.
The exception to this is the basis function scaling approach applied in the context of XFEM, which, in fact, can be interpreted as a diagonal preconditioning technique.
It is the primary objective of this work to extend this approach to be effective in the context of the (I)FCM,
yielding a fully automated and robust algebraic preconditioning technique (\emph{i.e.,} a preconditioner constructed exclusively with information from the system matrix).

In this work, we perform a detailed analysis of the conditioning problems associated with (I)FCM.
An important novel contribution of this work is the derivation of an explicit scaling relation between the condition number and the smallest cell volume fraction
for the discretization of elliptic second order partial differential equations using polynomial bases on uniform meshes.
This relation reveals a strong dependence of the condition number on the order of the employed discretization, which corroborates the need to develop a strategy to improve the conditioning of (I)FCM. Motivated by the aforementioned relation, in this work an algebraic preconditioning technique called SIPIC (\emph{Symmetric Incomplete Permuted Inverse Cholesky}) is developed, based on basis function scaling in combination with local orthonormalization on cells with very small volume fractions. This local orthonormalization approach is similar to the manipulation of basis functions in the Stable Generalized Finite Element Method \cite{Babuska2012}. We present the algorithm for the construction of an algebraic preconditioner, which generally preserves the sparsity pattern of the system or yields negligible fill-in.
The construction of this preconditioner does not add any significant computational cost. 

Section~\ref{sec:FCM} of this paper presents the variational formulation of the finite cell method.
In Section~\ref{sec:conditioning}, the explicit relation between the condition number and the smallest cell volume fraction for uniform meshes is derived and numerically verified.
The construction of the algebraic preconditioner is described in Section~\ref{sec:precon}.
In Section~\ref{sec:numerical}, the numerical implementation is demonstrated and the effect of the proposed preconditioner is shown for numerical examples.
A novel contribution in this section is the weak imposition of Dirichlet boundary conditions in either only the normal or only the tangential direction.
Furthermore, it is demonstrated that preconditioning enables detailed mesh convergence analysis of FCM.
Conclusions are finally drawn in Section~\ref{sec:conclusion}.
A special conditioning technique for the local eigenvalue problem is described in \ref{sec:local}.

\section{The finite cell method}
\label{sec:FCM}

\subsection{Model problem and weak form}\label{sec:variational}
\begin{figure}
  \begin{center}
  \includegraphics[width=0.65\textwidth]{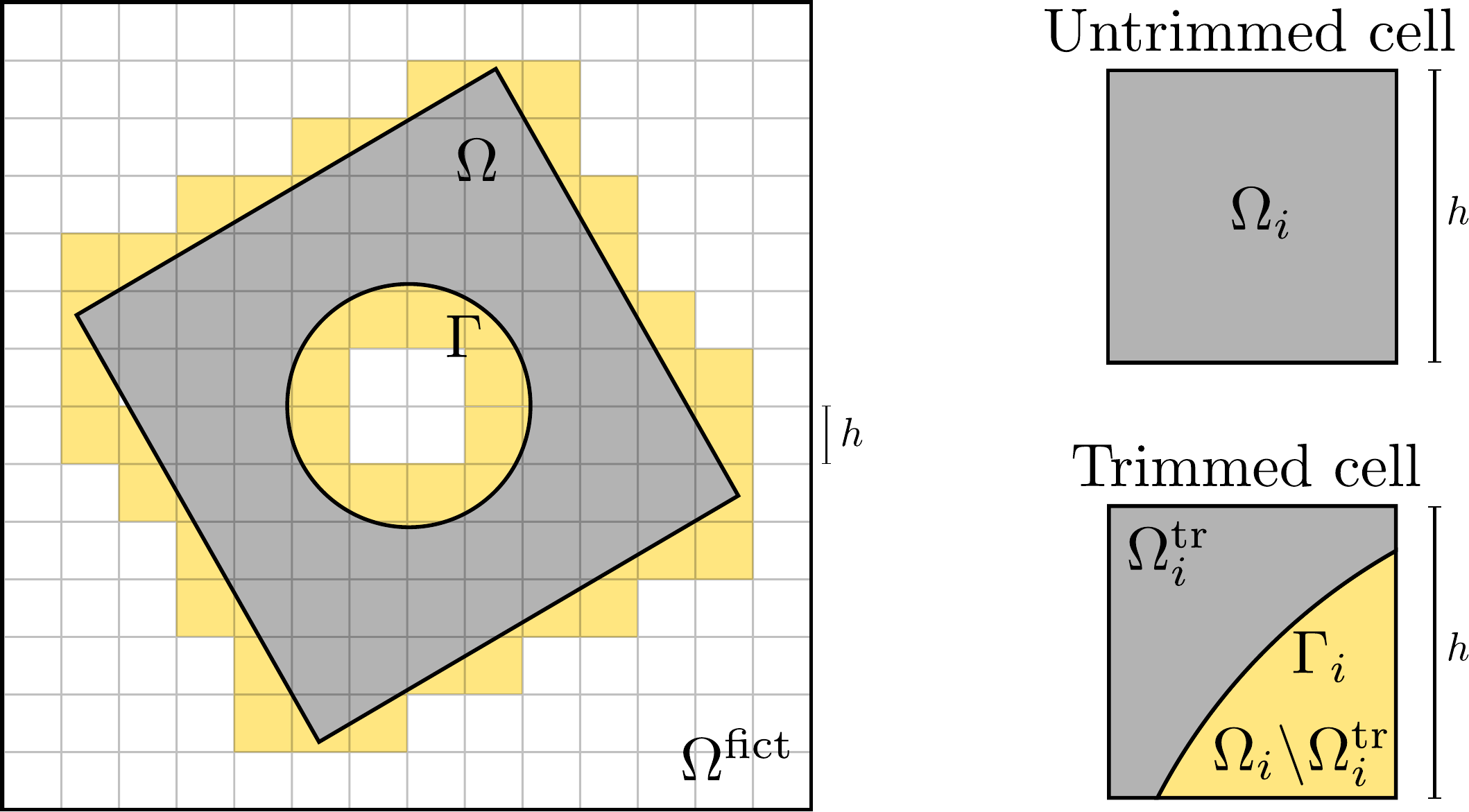}
  \caption{A geometrically complex domain $\Omega$ that is encapsulated by the geometrically simple, rectilinear domain $\Omega\cup\Omega_{\rm fict}$.\label{fig:domain}}
  \end{center}
\end{figure}

As a model problem we consider Poisson's equation over an open bounded Lipschitz domain $\Omega \subset \mathbb{R}^d$ ($d\in\{2,3\}$), supplemented with Dirichlet and Neumann conditions on complementary parts of the boundary $\partial \Omega$:
\begin{equation}
 \begin{cases}
  -\Delta u = f & {\rm in \ } \Omega, \\
  \partial_n u = g^N & {\rm on \ } \Gamma^N, \\
  u = g^D & {\rm on \ } \Gamma^D.     
 \end{cases}
 \label{eq:problem}
\end{equation}
In this expression, $\partial_n=n\cdot \nabla$ denotes the gradient in the direction of the the outward pointing normal vector $n$, and furthermore $\Gamma^D \cap \Gamma^N = \emptyset$ and $\overline{\Gamma^N \cup \Gamma^D} = \Gamma = \partial\Omega$.
Classical Galerkin methods for solving this problem (such as the Finite Element Method (FEM), \emph{e.g.,} \citep{FEM_hughes,FEM_zienkiewics}) employ boundary fitted meshes on which essential boundary conditions can be imposed in a strong manner -- \emph{i.e.,} encoded into the approximation space.
The Finite Cell Method (FCM) is an unfitted Galerkin method, and uses a geometrically simple mesh that encapsulates $\Omega$ such as a rectilinear discretization of $\Omega\cup\Omega_{\rm fict}$ as shown in Figure~\ref{fig:domain}.

To derive a consistent variational form for \eqref{eq:problem}, we start with the weighted residual form with test function $v \in H^1(\Omega)$\footnote{$H^1(\Omega)$ denotes the usual set of all square integrable functions with square integrable derivatives.}:
\begin{equation}
\int_\Omega -v \Delta u {\rm d}V = \int_\Omega v f {\rm d}V.
\label{eq:weighted_residual}
\end{equation}
Assuming that $\partial_n u$ exists on $\Gamma^D$, integration-by-parts in combination with the Neumann condition in \eqref{eq:problem} yields
\begin{equation}
\int_\Omega \nabla v \cdot \nabla u {\rm d}V + \int_{\Gamma^D} - v \partial_n u {\rm d}S = \int_\Omega v f {\rm d}V + \int_{\Gamma^N} vg^N {\rm d}S.
\label{eq:bulk}
\end{equation}
In most classical weak forms the integral over $\Gamma^D$ is left out, because in these formulations $v$ is chosen from a subset of $H^1(\Omega)$ such that $v=0$ on $\Gamma^D$.
This is not straightforward on unfitted meshes however.
Under the assumption that $\partial_n v$ exists on $\Gamma^D$, one can add the term 
\begin{equation}
\int_{\Gamma^D} - u \partial_n v {\rm d}S = \int_{\Gamma^D} - g^D \partial_n v {\rm d}S,
\label{eq:boundary}
\end{equation}
which is consistent with \eqref{eq:problem} on account of the essential boundary condition on $\Gamma^{D}$. The modification in \eqref{eq:boundary} leads to symmetry of the variational form \citep{Nitsche}:
\begin{equation}\begin{aligned}
& \int_\Omega \nabla v \cdot \nabla u {\rm d}V + \int_{\Gamma^D} - \left( v \partial_n u + u \partial_n v\right) {\rm d}S \\ \hspace{1cm} &= \int_\Omega v f {\rm d}V + \int_{\Gamma^N} vg^N {\rm d}S + \int_{\Gamma^D} - g^D \partial_n v {\rm d}S.
\end{aligned}\label{eq:consistent}
\end{equation}
Note that $\partial_n v$ exists on $\Gamma^D$ for all element-wise polynomial approximation spaces, but that this formulation is unbounded in the full space $H^1(\Omega)$. Because the boundary term in \eqref{eq:consistent} is not necessarily positive, we add a penalty term to ensure coercivity of the variational form (and consequently positive definiteness of the system matrix):
\begin{equation}
\int_{\Gamma^D} \beta vu {\rm d}S = \int_{\Gamma^D} \beta vg^D {\rm d}S.
\label{eq:penalty}
\end{equation}
The required value for the penalty parameter $\beta$ depends on the employed approximation space and -- because of the unbounded nature of the boundary term in \eqref{eq:consistent} on $H^1(\Omega)$ -- there exists no finite $\beta$ which is uniformly valid for all approximation spaces.

Considering a finite dimensional approximation space $\mathcal{V}_h(\Omega) \subset H^1(\Omega)$, the weak formulation of the approximation to \eqref{eq:problem}, based on \eqref{eq:weighted_residual}--\eqref{eq:penalty}, writes
\begin{equation}\left\{ \begin{array}{l}
\mbox{find } u_h \in \mathcal{V}_h(\Omega) \mbox{ such that:} \\[0.1em]
\mathcal{F}(v_h,u_h) = \ell (v_h) \quad \forall v_h \in \mathcal{V}_h(\Omega)
\end{array}\right.
\label{eq:weakformproblem}
\end{equation}
with
\begin{subequations}\label{eq:bilop}\begin{align}
\mathcal{F}(v_h,u_h) & = \int_\Omega \nabla v_h \cdot \nabla u_h {\rm d}V + \int_{\Gamma^D} \big( \beta v_hu_h - \left( v_h \partial_n u_h + u_h \partial_n v_h \right) \big) {\rm d}S, \\
\ell(v_h)   & = \int_\Omega v_h f {\rm d}V + \int_{\Gamma^N} v_hg^N {\rm d}S + \int_{\Gamma^D} \big( \beta v_hg^D - g^D \partial_n v_h \big) {\rm d}S.
\end{align}\end{subequations}
For notational convenience, the bilinear operator $\mathcal{F}(\cdot,\cdot)$ is split into three components,\\
$\mathcal{F}(\cdot,\cdot) = \mathcal{F}^1(\cdot,\cdot) + \mathcal{F}^2(\cdot,\cdot) + \mathcal{F}^3(\cdot,\cdot)$, with
\begin{subequations}\begin{align}
\mathcal{F}^1(v_h,u_h) & = \int_\Omega \nabla v_h \cdot \nabla u_h {\rm d}V, \\
\mathcal{F}^2(v_h,u_h) & = \int_{\Gamma^D} -\left( v_h \partial_n u_h + u_h \partial_n v_h \right) {\rm d}S, \\
\mathcal{F}^3(v_h,u_h) & = \int_{\Gamma^D} \beta v_hu_h {\rm d}S.
\end{align}\end{subequations}
Because $\mathcal{F}(\cdot,\cdot)$ is coercive on finite dimensional function spaces with respect to the $H^1(\Omega)$ norm, it defines the equivalent norm
\begin{equation}
\|v_h\|_{\mathcal{F}}^2 = \mathcal{F}(v_h,v_h),
\end{equation}
referred to as the finite cell norm (or FCM norm).

\subsection{Continuity and coercivity of the discrete bilinear form}
\label{sec:coercivity}
The stability parameter $\beta$ -- used to ensure coercivity of the variational form, see \eqref{eq:penalty} -- can be chosen globally or locally (\emph{i.e.,} one constant for the whole boundary $\Gamma^D$ or a separate constant for every cell that is intersected by $\Gamma^D$). The global approach is presented in \citep{Embar2010}, and is identical to the local approach when the whole domain is treated as one cell. A standard discrete trace inequality (\emph{e.g.,} \cite{InverseTrace}) conveys that with the local approach, the lower bound for $\beta_i$ on cell $\Omega_i$ is of the order of magnitude of $1/\widetilde{h}_i$, with $\widetilde{h}_i$ the typical length scale of $\Omega_i^{\rm tr}=\Omega_i\cap\Omega$ (Figure~\ref{fig:domain}). With the global approach, the lower bound for $\beta$ is of the same order of magnitude as $\max_i \beta_i$. Consequently, one cell that is intersected by $\Gamma^D$ such that its $\Omega_i^{\rm tr}$ is small, will cause a large value of $\beta$ on the entire boundary $\Gamma^D$. As a result, global stabilization negatively affects the conditioning of FCM (\ref{sec:app_bound}). Moreover, it can degenerate FCM to a penalty method \citep{Ruess2013,Ruess2014}. Local stabilization -- as employed herein -- is not prone to these deficiencies.

For $\mathcal{F}(\cdot,\cdot)$ to be coercive, on every cell it must hold that $\beta_i > C_i$, with \citep{Embar2010}:
\begin{equation}
C_i = \max_{v_h \in \mathcal{V}_{h}|_{\Omega_i}} \frac{\| \partial_n v_h \|_{L^2(\Gamma_i^D)}^2}{\mathcal{F}^1_i(v_h,v_h)},
\label{eq:infdim}\end{equation}
where $\mathcal{V}_{h}|_{\Omega_i}$ denotes the restriction of the approximation space $\mathcal{V}_h$ to $\Omega_i$, $\Gamma_i^D = \Omega_i\cap\Gamma^D$ is the part part of the boundary $\Gamma^D$ contained in cell $\Omega_i$ (Figure~\ref{fig:domain}), and
\begin{equation}
\mathcal{F}^1_i(v_h,u_h) =   \int_{\Omega_i^{\rm tr}} \nabla v_h \cdot \nabla u_h {\rm d}V.
\end{equation}
Coercivity of $\mathcal{F}_i(\cdot,\cdot)$ subject to the condition $\beta_i > C_i$ follows from the following sequence of inequalities:
\begin{equation}
\begin{aligned}
\mathcal{F}_i(v_h,v_h) = & \mathcal{F}^1_i(v_h,v_h)+\mathcal{F}^2_i(v_h,v_h)+\mathcal{F}^3_i(v_h,v_h) \\
\geq & \mathcal{F}^1_i(v_h,v_h)+\mathcal{F}^3_i(v_h,v_h)-|\mathcal{F}^2_i(v_h,v_h)| \\
\geq & \mathcal{F}^1_i(v_h,v_h)+\mathcal{F}^3_i(v_h,v_h) - 2\| n \cdot \nabla v_h \|_{L^2(\Gamma_i^D)}\| v_h \|_{L^2(\Gamma_i^D)}\\
\geq & \mathcal{F}^1_i(v_h,v_h)+\mathcal{F}^3_i(v_h,v_h) - \varepsilon_i \| n \cdot \nabla v_h \|_{L^2(\Gamma_i^D)}^2 - {\varepsilon_i^{-1}} \| v_h \|_{L^2(\Gamma_i^D)}^2\\
\geq & \mathcal{F}^1_i(v_h,v_h)+\mathcal{F}^3_i(v_h,v_h) - C_i\varepsilon_i \mathcal{F}^1_i(v_h,v_h) + {\varepsilon_i^{-1}} \| v_h \|_{L^2(\Gamma_i^D)}^2\\
\geq & \left( 1 - C_i\varepsilon_i \right) \mathcal{F}^1_i(v_h,v_h) + \int_{\Gamma_i^D} (\beta_i-{\varepsilon_i^{-1}}) v_h^2 {\rm d}S\\
\geq & \delta_i \|v_h\|_{H^1(\Omega_i^{\rm tr})}^2 > 0,
\end{aligned}
\label{eq:proof_coercivity}
\end{equation}
for some $\delta_i>0$. The third step in \eqref{eq:proof_coercivity} follows from the Cauchy-Schwarz inequality. Step four follows from Young's inequality, while  equation~\eqref{eq:infdim} is substituted in step five. Step six is based on the fact that there exists a constant $\varepsilon_i$ such that $C_i < \varepsilon_i^{-1} < \beta_i$. The last step applies a specific form of the Poincar\'e inequality (\emph{e.g.,} Lemma B.63 in \citep{ErnGuermond}). Global coercivity, $\mathcal{F}(v_h,v_h) \geq \delta \|v_h\|_{H^1(\Omega)}^2 > 0$ for some $\delta>0$, follows by summation over all $\mathcal{F}_i$. In this work, $\beta_i = 2 C_i$ is applied for all examples and test cases. This choice is essentially arbitrary and does not affect the strategy to construct a preconditioner.

The values of $C_i$ can be computed numerically. Because $\mathcal{V}_{h}|_{\Omega_i}$ is finite dimensional, \eqref{eq:infdim} can be cast in matrix form as
\begin{equation}
C_i = \max_{\mathbf{x}} \frac{\mathbf{x}^T \mathbf{B}_i \mathbf{x}}{\mathbf{x}^T \mathbf{V}_i \mathbf{x}},
\label{eq:findim}\end{equation}
with
\begin{subequations}\begin{align}
\mathbf{B}_i &= \int_{\Gamma_i^D} \partial_n \boldsymbol{\Phi}_i \partial_n \boldsymbol{\Phi}_i^T {\rm d}S, \label{eq:localevB}\\
\mathbf{V}_i &= \mathcal{F}^1_i(\boldsymbol{\Phi}_i,\boldsymbol{\Phi}_i^T) = \int_{\Omega_i^{\rm tr}} \nabla \boldsymbol{\Phi}_i \cdot \nabla \boldsymbol{\Phi}_i^T {\rm d}V,\label{eq:localevV}
\end{align}\end{subequations}
where $\boldsymbol{\Phi}_i$ denotes the vector of all basis functions that are supported on $\Omega_i$.
The optimality condition associated with \eqref{eq:findim} is
\begin{equation}
\partial_{\mathbf{x}} \frac{\mathbf{x}^T \mathbf{B}_i \mathbf{x}}{\mathbf{x}^T \mathbf{V}_i \mathbf{x}} = \frac{\left(\mathbf{B}_i + \mathbf{B}_i^T\right) - \frac{\mathbf{x}^T \mathbf{B}_i \mathbf{x}}{\mathbf{x}^T \mathbf{V}_i \mathbf{x}}\left(\mathbf{V}_i + \mathbf{V}_i^T\right) }{\mathbf{x}^T\mathbf{V}_i\mathbf{x}} \mathbf{x} = \mathbf{0},
\end{equation}
which implies
\begin{equation}
\left(\mathbf{B}_i - \frac{\mathbf{x}^T \mathbf{B}_i \mathbf{x}}{\mathbf{x}^T \mathbf{V}_i \mathbf{x}} \mathbf{V}_i \right) \mathbf{x} = \mathbf{0},
\end{equation}
because the matrices $\mathbf{B}_i$ and $\mathbf{V}_i$ are both symmetric.
The solution to \eqref{eq:infdim} therefore coincides with the largest eigenvalue of the generalized eigenvalue problem
\begin{equation}
\mathbf{B}_i\mathbf{x} = \lambda \mathbf{V}_i\mathbf{x},
\label{eq:geneig}
\end{equation}
which enables computation of $C_i=\lambda_{\rm max}$.

The generalized eigenvalue problem \eqref{eq:geneig} may suffer from conditioning problems. A solution method that ameliorates these conditioning problems is discussed in \ref{sec:local}.

\section{Condition number analysis}
\label{sec:conditioning}
The condition number of an invertible matrix $\mathbf{A}$ is defined as
\begin{equation}
\kappa_2(\mathbf{A})= \|\mathbf{A}\|_2 \|\mathbf{A}^{-1}\|_2 \geq 1,
\end{equation}
with $\|\cdot\|_2$ the (induced) Euclidean norm:
\begin{align}
\|\mathbf{A}\|_2 = \max_{\mathbf{x}} \frac{\|\mathbf{A}\mathbf{x}\|_2}{\|\mathbf{x}\|_2}.
\label{eq:ind_norm}
\end{align}
When solving a linear system of the form
\begin{equation}
\mathbf{A}\mathbf{x} = \mathbf{b},
\end{equation}
the condition number $\kappa_2(\mathbf{A})$ is of major importance. 
Well-conditioned systems (condition numbers close to $1$) have a smaller propagation of errors in $\mathbf{A}$ or $\mathbf{b}$,
converge faster when solved iteratively, and have a smaller uncertainty when the error is estimated using the residual.
For example, when $\mathbf{A}$ is symmetric positive definite -- such as system matrices originating from most FCM approximations -- 
the Conjugate Gradient method has a convergence bound that depends on the condition number as
\begin{equation}
\|\mathbf{x} - \mathbf{x}_i\|_\mathbf{A} \leq 2 \left(\frac{\sqrt{\kappa_2(\mathbf{A})}-1}{\sqrt{\kappa_2(\mathbf{A})}+1}\right)^i \|\mathbf{x} - \mathbf{x}_0\|_\mathbf{A},
\label{eq:cgconvergence}
\end{equation}
where $\mathbf{x}_i$ is the approximation after $i$ iterations and $\|\cdot\|_\mathbf{A}$ denotes the energy norm $\|\mathbf{x}\|_\mathbf{A}^2 = \mathbf{x}^T \mathbf{A}\mathbf{x}$ \citep{Saad}.
Hence, the guaranteed reduction in the error per iteration is much larger for $\kappa_2(\mathbf{A}) \sim 1$ than for $\kappa_2(\mathbf{A}) \gg 1$.
The estimate of $\|\mathbf{x} - \mathbf{x}_i\|_\mathbf{A}$ in terms of the residual $\| \mathbf{b} - \mathbf{A} \mathbf{x}_i \|_2$ also depends on the condition number:
\begin{equation}
\frac{\| \mathbf{b} - \mathbf{A} \mathbf{x}_i \|_2^2}{\|\mathbf{A}\|_2}  \leq \| \mathbf{x} - \mathbf{x}_i \|_\mathbf{A}^2 \leq \kappa_2(\mathbf{A})\frac{\| \mathbf{b} - \mathbf{A} \mathbf{x}_i \|_2^2}{\| \mathbf{A} \|_2}.
\label{eq:cgerrorbounds}
\end{equation}
Equation~\eqref{eq:cgerrorbounds} conveys that, in the case of ill-conditioning, reduction of the residual to a value below a specified tolerance does not guarantee convergence by the same amount of the (unknown) error in the solution (\emph{e.g.,} the error in the energy norm).

Because the bilinear operator in \eqref{eq:bilop} is symmetric and coercive, system matrices originating from this FCM formulation are Symmetric Positive Definite (SPD). In Section~\ref{sec:cond_SPD} we derive an equivalent definition of the norm for SPD matrices, which enables the application to FCM system matrices in Section~\ref{sec:cond_FCM}. In Section~\ref{sec:cond_example} we verify the obtained theoretical estimate of the condition number with an example.

\subsection{Equivalent norms for SPD matrices}\label{sec:cond_SPD}
Because symmetric matrices have orthogonal eigenvectors, definition~\eqref{eq:ind_norm} implies
\begin{equation}
\|\mathbf{A}\|_2 = |\lambda|_{\rm max} \qquad\mbox{and}\qquad \|\mathbf{A}^{-1}\|_2 = \frac{1}{|\lambda|_{\rm min}},
\label{eq:norm1}\end{equation}
in which $|\lambda|_{\rm max}$ and $|\lambda|_{\rm min}$ denote the absolute maximal and minimal eigenvalues of $\mathbf{A}$, respectively.
Furthermore, for symmetric matrices the Rayleigh quotient is bounded by the eigenvalues,
\begin{equation}
 \lambda_{\rm min} = \frac{\mathbf{y}^T\lambda_{\rm min}\mathbf{y}}{\mathbf{y}^T\mathbf{y}} \leq
 \frac{\mathbf{y}^T\mathbf{A}\mathbf{y}}{\mathbf{y}^T\mathbf{y}} \leq \frac{\mathbf{y}^T\lambda_{\rm max}\mathbf{y}}{\mathbf{y}^T\mathbf{y}} = \lambda_{\rm max},
\label{eq:norm2}\end{equation}
where the left, respectively right, inequality becomes sharp when $\mathbf{y}$ coincides with the eigenvector corresponding to the smallest, respectively largest, eigenvalue.
As a result,
\begin{equation}
 \lambda_{\rm max} = \max_{\mathbf{y}} \frac{\mathbf{y}^T\mathbf{A}\mathbf{y}}{\mathbf{y}^T\mathbf{y}} \qquad\mbox{and}\qquad \lambda_{\rm min} =
 \min_{\mathbf{y}} \frac{\mathbf{y}^T\mathbf{A}\mathbf{y}}{\mathbf{y}^T\mathbf{y}}.
\label{eq:norm3}\end{equation}
If $\mathbf{A}$ is positive definite as well (\emph{i.e.,} $\mathbf{A}$ is SPD), all eigenvalues are positive such that
\begin{equation}
\lambda_{\rm max} = |\lambda|_{\rm max} \qquad\mbox{and}\qquad \lambda_{\rm min} = |\lambda|_{\rm min}. 
\label{eq:norm4}\end{equation}
Combination of \eqref{eq:norm1}--\eqref{eq:norm4} yields
\begin{equation}
 \|\mathbf{A}\|_2 = \max_{\mathbf{y}} \frac{\mathbf{y}^T\mathbf{A}\mathbf{y}}{\mathbf{y}^T\mathbf{y}} \qquad\mbox{and}\qquad \|\mathbf{A}^{-1}\|_2 =
 \max_{\mathbf{y}} \frac{\mathbf{y}^T\mathbf{y}}{\mathbf{y}^T\mathbf{A}\mathbf{y}},
\label{eq:normratio}\end{equation}
which shows that the norm and inverse norm of SPD matrices are determined by the maximal and minimal quotient of the energy norm and the Euclidean norm of a vector.

\subsection{Condition numbers in FCM}\label{sec:cond_FCM}
To apply \eqref{eq:normratio} to system matrices originating from FCM formulations similar to \eqref{eq:bilop},
we consider a function $v_h = \boldsymbol{\Phi}^T \mathbf{y} \in \mathcal{V}_h$ for some unique $\mathbf{y} \in \mathbb{R}^n$, with $\boldsymbol{\Phi}$ denoting the vector containing all basis functions.
The FCM system matrix is defined as
\begin{equation}
 \mathbf{A} = \mathcal{F}(\boldsymbol{\Phi},\boldsymbol{\Phi}^T),
\end{equation}
and the following is noted: 
\begin{equation}
\|\mathbf{y}\|_\mathbf{A}^2 = \mathbf{y}^T\mathbf{A}\mathbf{y} = \mathcal{F}(\mathbf{y}^T\boldsymbol{\Phi},\boldsymbol{\Phi}^T\mathbf{y}) = \mathcal{F}(v_h,v_h) = \|v_h\|_{\mathcal{F}}^2,
\end{equation}
showing that the FCM norm of a function $v_h$ is equal to the energy norm of the corresponding vector $\mathbf{y}$ with system matrix $\mathbf{A}$.
Therefore the $\| \cdot \|_2$ norm of $\mathbf{A}$ and its inverse according to \eqref{eq:normratio} can be interpreted as the quotient of the FCM norm of a function and the Euclidean norm of the corresponding vector:
\begin{equation}
  \|\mathbf{A}\|_2 = \max_{v_h,\mathbf{y}} \frac{\|v_h\|_{\mathcal{F}}^2}{\|\mathbf{y}\|_2^2}, \quad \|\mathbf{A}^{-1}\|_2 =
 \max_{v_h,\mathbf{y}} \frac{\|\mathbf{y}\|_2^2}{\|v_h\|_{\mathcal{F}}^2}.
\label{eq:funcvecratio}
\end{equation}
When a function $v_h$ exists for which the FCM norm $\|v_h\|_{\mathcal{F}}$ is very small in relation to the corresponding norm $\|\mathbf{y}\|_2$ of the coefficients vector, ill-conditioning can occur on account of the norm $\|\mathbf{A}^{-1}\|_2$ being very large. This is a situation typically encountered in FCM. Since in FCM the relative position of the mesh to the physical domain $\Omega$ is arbitrary, the volume of a trimmed cell $\Omega_i^{\rm tr}=\Omega_i\cap\Omega$ (Figure~\ref{fig:domain}) can be arbitrarily small. When a function is only supported on $\Omega_i^{\rm tr}$, its support\footnotemark -- and as a result also its FCM norm -- can therefore be arbitrarily small as well. \footnotetext{In the remainder of this manuscript, the term support refers to the support in the physical domain.} The vector that corresponds to this function does not depend on the volume of $\Omega_i^{\rm tr}$, however. As a result, the inverse norm $\|\mathbf{A}^{-1}\|_2$ in \eqref{eq:funcvecratio} can become arbitrarily large. In contrast, the norm $\|\mathbf{A}\|_2$ in itself is not sensitive to the relative position of the mesh. Consequently, FCM can yield arbitrarily large condition numbers, in case the mesh is positioned such that a, or multiple, small trimmed cells occur.

To provide a quantitative estimate of the condition number, we restrict ourselves to situations where:
\begin{itemize}
\item The number of dimensions $d$, is larger than one.
\item A local stabilization parameter $\beta_i$, is used, as described in Section~\ref{sec:coercivity}.
\item A piecewise polynomial basis is applied (\emph{e.g.,} B-splines, Lagrange, or spectral bases).
\item A uniform mesh with mesh size $h$ is used. 
\end{itemize}
The volume fraction $\eta_i$ is defined as the fraction of cell $\Omega_i$ that intersects physical domain $\Omega$,
\begin{equation}
\eta_i = \frac{|\Omega_i\cap\Omega|}{|\Omega_i|} = \frac{|\Omega_i^{\rm tr}|}{|\Omega_i|} = \frac{|\Omega_i^{\rm tr}|}{h^d}.
\label{eq:etai}
\end{equation}
The index of the smallest volume fraction is denoted by
\begin{equation}
\imath = \underset{i}{\operatorname{arg\,min}} ~\eta_i,
\end{equation}
such that the smallest volume fraction is $\eta = \eta_{\imath}$. Furthermore, we make three assumptions on the shape of the trimmed cells:
\begin{enumerate}
\item $\partial \Omega_{\imath} \cap \Omega$, the intersection between the untrimmed boundary of the cell with the smallest volume fraction and the physical domain $\Omega$,
contains parts of at most one of two opposing cell faces.
For example, if $\Omega_{\imath} = (0,1)^d$, then 
\begin{equation}
\left\{x \in \Omega_{\imath}^{\rm tr} | x_j = 0\right\} \neq \emptyset \quad \Rightarrow \quad \left\{x \in \Omega_{\imath}^{\rm tr} | x_j = 1\right\} = \emptyset,
\end{equation}
for all indices $j \leq d$.
\item There exists a positive constant $C_R$ such that for every cell $\Omega_i$, the radius $R_i$ of the smallest ball enclosing $\Omega_i^{\rm tr}$ is bounded by
\begin{equation}
R_i \leq C_R |\Omega_i^{\rm tr}|^{\frac{1}{d}}.
\end{equation}
It follows from \eqref{eq:etai} that $R_i \leq C_R h \eta_i^{\frac{1}{d}}$.
\item There exists a positive constant $C_\Gamma$ such that for every cell $\Omega_i$ that is trimmed by $\Gamma$, the surface measure of $\Gamma_i = \Gamma \cap \Omega_i$ (Figure~\ref{fig:domain}) is bounded by
\begin{equation}
|\Gamma_i| \leq C_\Gamma |\Omega_i^{\rm tr}|^{\frac{d-1}{d}} = C_\Gamma h^{d-1} \eta^{\frac{d-1}{d}}.
\end{equation}
\end{enumerate}
The first assumption implies that there exists a function in the approximation space that is only supported on $\Omega_{\imath}$. For $\eta$ small enough, this first assumption is automatically satisfied if assumption 2 or 3 holds. The second assumption guarantees shape regularity of trimmed cells, \emph{e.g.,} \citep{Hansbo2002,Burman2014}. This assumption underlies the FCM method, and serves to keep the penalty parameter bounded. This excludes pathological cases, such as trimmed cells with highly distorted aspect ratios. The third assumption bounds the size of the intersection between $\Gamma$ and a cell by the size of the trimmed cell. Under assumptions 2 and 3, it can be shown that $\exists C_\beta > 0$ such that the local stabilization parameter satisfies
\begin{equation}
\beta_i \leq C_\beta |\Omega_i^{\rm tr}|^{-1/d} = C_\beta h^{-1} \eta_i^{-1/d},
\label{eq:betabound}\end{equation}
see \emph{e.g.,} \citep{InverseTrace}. The restrictions and assumptions mentioned above, are only required to prove an estimate of the condition number of the FCM system matrix. The strategy to construct a preconditioner (proposed in Section~\ref{sec:precon}) does not depend on these restrictions and assumptions.

From the literature (\emph{e.g.,} \citep{johnson}) we know that the condition number of system matrices originating from classical FEM approximations of the Laplace operator with quasi-uniform meshes, scales with $h^{-2}$. The norm and largest eigenvalue of such a matrix scales with $h^{d-2}$ and the corresponding eigenvector represents the function with the highest frequency that can be represented on the mesh. The smallest eigenvalue of such a matrix scales with $h^d$ and the corresponding eigenvector represents the function with the lowest frequency. Note that this function approximates the analytical lowest eigenmode of the operator, which is independent of the mesh size $h$. However, the norm of the corresponding coefficient vector does depend on $h$ by virtue of the fact that the dimension of the approximation space is altered under mesh refinement. Because $\mathcal{F}^1(\cdot,\cdot)$ coincides with the standard FEM bilinear form, for the norm $\|\mathbf{A}\|_2$ it then follows that\footnote{In the remainder of this manuscript, the variable $C$, without subscript or superscript, denotes a general positive constant that may attain different values in different statements or equations.}
\begin{equation}\begin{aligned}
 \|\mathbf{A}\|_2 & = \max_{\mathbf{y}} \frac{\mathbf{y}^T\mathbf{A}\mathbf{y}}{\mathbf{y}^T\mathbf{y}} =
 \max_{v_h,\mathbf{y}} \frac{\mathcal{F}(v_h,v_h)}{\mathbf{y}^T\mathbf{y}} \geq \max_{v_h|_{\Gamma^D}=0,\mathbf{y}} \frac{\mathcal{F}(v_h,v_h)}{\mathbf{y}^T\mathbf{y}} \\ & = \max_{v_h|_{\Gamma^D}=0,\mathbf{y}} \frac{\mathcal{F}^1(v_h,v_h)}{\mathbf{y}^T\mathbf{y}} \geq C h^{d-2}.
 \label{eq:femeigenvalue}
\end{aligned}\end{equation}
Under assumptions 2 and 3, it can be verified that (see \ref{sec:app_bound} for details)
\begin{align}
&\frac{|\mathcal{F}^2(v_h,v_h)|}{\mathbf{y}^T\mathbf{y}}  \leq C h^{d-2}  \qquad\mbox{and}\qquad \frac{\mathcal{F}^3(v_h,v_h)}{\mathbf{y}^T\mathbf{y}}  \leq C h^{d-2}.
\label{eq:bcbounds}
\end{align}
From \eqref{eq:femeigenvalue} and \eqref{eq:bcbounds} we can derive the upper bound
\begin{equation}
\|\mathbf{A}\|_2 \geq c_\mathcal{F} h^{d-2},
\label{eq:est_norm}
\end{equation}
for some $c_\mathcal{F}>0$. We emphasize that \eqref{eq:bcbounds} only holds for locally stabilized systems. A derivation of \eqref{eq:bcbounds} and an elaboration of the effect of global stabilization is given in \ref{sec:app_bound}.

To estimate $\|\mathbf{A}^{-1}\|$ using \eqref{eq:funcvecratio}, we need to evaluate the smallest eigenvalue.
As mentioned before, for matrices originating from classical FEM formulations, the eigenvector that corresponds to the smallest eigenvalue represents the function with the lowest frequency.
This is generally not the case for matrices originating from FCM formulations.
In most FCM cases, the smallest eigenvalue corresponds to a function which is only supported on a cell with a very small volume fraction. 
Under assumption 1, there exists a function $v_h$ that is only supported on ${\Omega_{\imath}^{\rm tr}}$, 
\begin{equation}
v_h|_{\Omega_{\imath}^{\rm tr}} = \prod_{j=1}^d \left(\frac{x_j-\hat{x}_j}{h}\right)^{p},
\label{eq:smallfunc}
\end{equation}
which has a corresponding vector with a magnitude of order one, $\|\mathbf{y}\|_2 \sim 1$. In \eqref{eq:smallfunc}, $p$ is the order of the discretization, $j$ is the index of the dimension and $\hat{x} \in \partial\Omega_{\imath}\cap B_{R_{\imath}}(\Omega_{\imath}^{\rm tr})$
-- \emph{i.e.,} $\hat{x}$ lies in the intersection between the (untrimmed) boundary of $\Omega_{\imath}$ and the smallest ball enclosing $\Omega_{\imath}^{\rm tr}$.
For example, if there is a vertex in $\partial\Omega_{\imath}\cap B_{R_{\imath}}(\Omega_{\imath}^{\rm tr})$ (under assumption 1 there can be at most one), then $\hat{x}$ coincides with this vertex.
Under assumption 2, we can show that for the function $v_h$ according to \eqref{eq:smallfunc}, it holds that
\begin{equation}
\|v_h\|_{L^\infty\left(\Omega_{\imath}^{\rm tr}\right)} \leq \left(\frac{2R_{\imath}}{\sqrt{d}h}\right)^{pd} \leq \left(\frac{2C_R}{\sqrt{d}}\right)^{pd} \eta^{p},
\end{equation}
and
\begin{equation}
\| \nabla v_h \|_{L^\infty\left(\Omega_{\imath}^{\rm tr}\right)} \leq
\frac{p\sqrt{d}}{h} \left(\frac{2R_{\imath}}{\sqrt{d}h}\right)^{pd-1} \leq p\sqrt{d} \left(\frac{2C_R}{\sqrt{d}}\right)^{pd-1} \frac{\eta^{p-1/d}}{h},
\end{equation}
where $\| \nabla v_h \|_{L^\infty}$ denotes the supremum of the Euclidean norm of $\nabla v_h$. Recalling assumption 3 and the upper bound for $\beta_{\imath}$ given in \eqref{eq:betabound}, the following bounds hold
\begin{subequations}\begin{align}
 \mathcal{F}^1(v_h,v_h) & = \int_{\Omega_{\imath}^{\rm tr}} \nabla v_h \cdot \nabla v_h {\rm d}V \leq
 |\Omega_{\imath}^{\rm tr}| \| \nabla v_h \|_{L^\infty\left(\Omega_{\imath}^{\rm tr}\right)}^2 \\
 & \leq p^2d \left(\frac{2C_R}{\sqrt{d}}\right)^{2pd-2} h^{d-2} \eta^{2p+1-2/d}, \nonumber \\
 \mathcal{F}^2(v_h,v_h) & \leq \int_{\Gamma_{\imath}^D} 2 v_h |\partial_n v_h| {\rm d}S \\
 & \leq 2 |\Gamma_{\imath}^D| \|v_h\|_{L^\infty\left(\Omega_{\imath}^{\rm tr}\right)} \| \nabla v_h \|_{L^\infty\left(\Omega_{\imath}^{\rm tr}\right)} \nonumber \\
 & \leq 2 C_\Gamma p \sqrt{d} \left(\frac{2C_R}{\sqrt{d}}\right)^{2pd-1} h^{d-2} \eta^{2p+1-2/d}, \nonumber \\
 \mathcal{F}^3(v_h,v_h) & = \int_{\Gamma_{\imath}^D} \beta_{\imath} v_h^2 {\rm d}S \\
 & \leq |\Gamma_{\imath}^D| \beta_{\imath} \|v_h\|_{L^\infty\left(\Omega_{\imath}^{\rm tr}\right)}^2 \nonumber \\
 & \leq C_\Gamma C_\beta \left(\frac{2C_R}{\sqrt{d}}\right)^{2pd} h^{d-2} \eta^{2p+1-2/d}, \nonumber
\end{align}\end{subequations}
which enables us to give an upper bound for the complete FCM norm of the function specified in \eqref{eq:smallfunc}:
\begin{equation}
\|v_h\|_{\mathcal{F}}^2 \leq C_\mathcal{F} h^{d-2} \eta^{2p+1-2/d},
\end{equation}
with $C_\mathcal{F}$ according to
\begin{equation}
C_\mathcal{F} = p^2d \left(\frac{2C_R}{\sqrt{d}}\right)^{2pd-2} + 2 C_\Gamma p \sqrt{d} \left(\frac{2C_R}{\sqrt{d}}\right)^{2pd-1} + C_\Gamma C_\beta \left(\frac{2C_R}{\sqrt{d}}\right)^{2pd}.
\end{equation}
Because $\|\mathbf{y}\|_2 \sim 1$ for the function specified in \eqref{eq:smallfunc}, there exists a constant $c_{\mathcal{F}^{-1}} \sim 1/C_\mathcal{F}$ such that
\begin{equation}
\|\mathbf{A}^{-1}\| \geq c_{\mathcal{F}^{-1}} \frac{\eta^{-(2p+1-2/d)}}{h^{d-2}}.
\label{eq:est_invnorm}
\end{equation}
Combining \eqref{eq:est_norm} and \eqref{eq:est_invnorm}, we obtain the following \emph{lower bound} for the FCM condition number: 
\begin{equation}
\kappa_2(\mathbf{A}) \geq c_\mathcal{F} c_{\mathcal{F}^{-1}} \eta^{-(2p+1-2/d)}.
\label{eq:estimate}
\end{equation}
The condition number of matrices originating from FCM formulations is therefore bounded from below by a constant scaling with $\eta^{-(2p+1-2/d)}$. The dependence of the scaling rate on the order of the approximation space conveys that the combination of FCM with higher-order methods, such as p-FEM and IGA, makes it particularly susceptible to conditioning problems.

\subsection{Numerical test case: Unit square with circular exclusion}
\label{sec:cond_example}
To verify the scaling relation \eqref{eq:estimate}, we consider the finite cell discretization of problem \eqref{eq:problem} with Dirichlet boundary conditions on all boundaries ($\Gamma^D = \partial\Omega$ and $\Gamma^N = \emptyset$) on a unit square domain ($L=1$) with a centered circular exclusion of radius $R$, \emph{i.e.,} $\Omega = \{ x \in (-\frac{1}{2},\frac{1}{2})^2 :  | x | > R  \}$. Note that the data $f$ and $g_D$ in \eqref{eq:problem} do not appear in the system matrix and hence are not required for the computation of the condition number.

\begin{figure}
 \centering
 \includegraphics[width=0.5\linewidth]{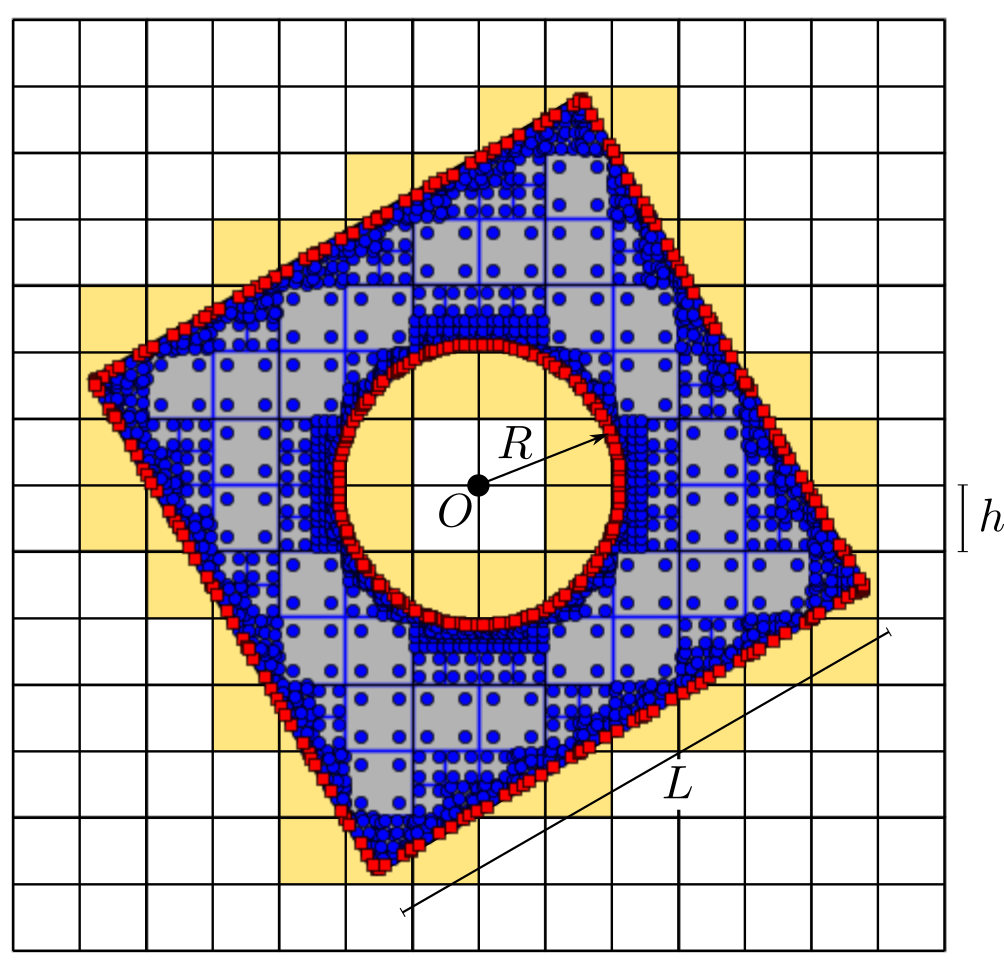}
 \caption{Schematic representation of a unit square domain with circular exclusion. Cells that intersect the physical domain are depicted in yellow. Volumetric (Gauss) integration points are indicated in blue and boundary integration points in red.}
 \label{fig:testcase_square}
\end{figure}

We partition the embedding domain ($\Omega \cup \Omega_{\rm fict}$) with a Cartesian mesh with mesh size $h=\frac{1}{32}$, with the center of the exclusion positioned at a vertex of the mesh. This is illustrated in Figure~\ref{fig:testcase_square}, where the cells in white do not intersect the physical domain and hence do not contribute to the system matrix. We consider uniform bivariate B-spline bases \cite{Hughes2005} of order $p\in\{1,2,3,4\}$ and Lagrange bases of order $p \in \{ 1, 2 \}$. The bisection-based tessellation scheme proposed in \citep{Verhoosel2015} with a maximal refinement depth of two is employed to accurately approximate the geometry of the domain. The number of Gauss points is selected such that exact integration of all operators is achieved over the tessellated domain. An applied illustration of the interior (blue circles) and boundary (red squares) integration points is shown in Figure~\ref{fig:testcase_square}. The results presented in this section pertain to finite cell systems with local stabilization parameter $\beta_i = 2 C_i$ (see Section~\ref{sec:coercivity}). Results for global stabilization -- although not presented here -- demonstrate the validity of scaling relation \eqref{eq:estimate_global}.

To construct scenarios with different volume fractions, the domain is gradually rotated about the origin while the background mesh remains unmodified. Starting from a domain in which the outer boundaries are aligned with the mesh, it is rotated in 100 steps over an angle of $\SI{45}{\degree}$. The radius $R$ is taken equal to $R=\sqrt{1/8-\sqrt{\eta_{R}/2}h} \approx \frac{1}{4}\sqrt{2}$, such that a volume fraction of approximately $\eta_R = 5 \cdot 10^{-3}$ occurs at the boundary of the circular exclusion, regardless of the angle over which the domain is rotated. The minimal volume fraction $\eta$, is therefore bounded from above by approximately $\eta_R$. The minimal volume fraction and condition number are computed for each configuration.

\begin{figure}
 \centering
 \begin{subfigure}{0.49\linewidth}
   \centering
   \includegraphics[width=\linewidth]{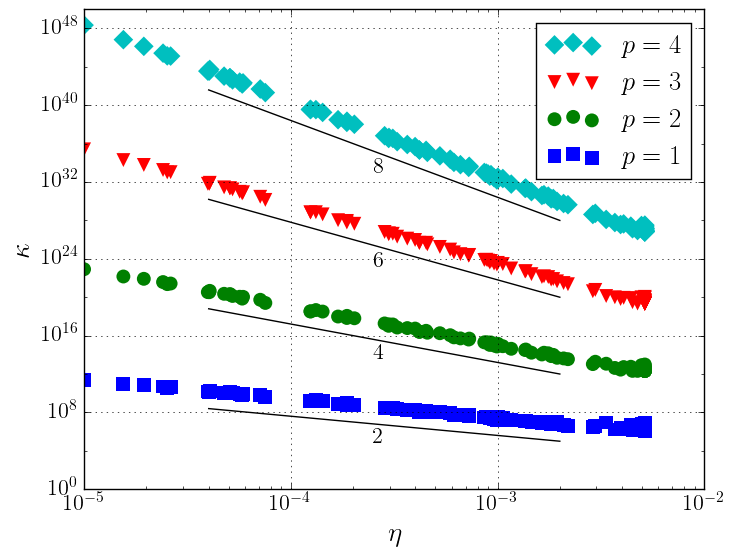}
   \caption{B-splines}
   \label{fig:condition_spline}
 \end{subfigure}
 \begin{subfigure}{0.49\linewidth}
   \centering
   \includegraphics[width=\linewidth]{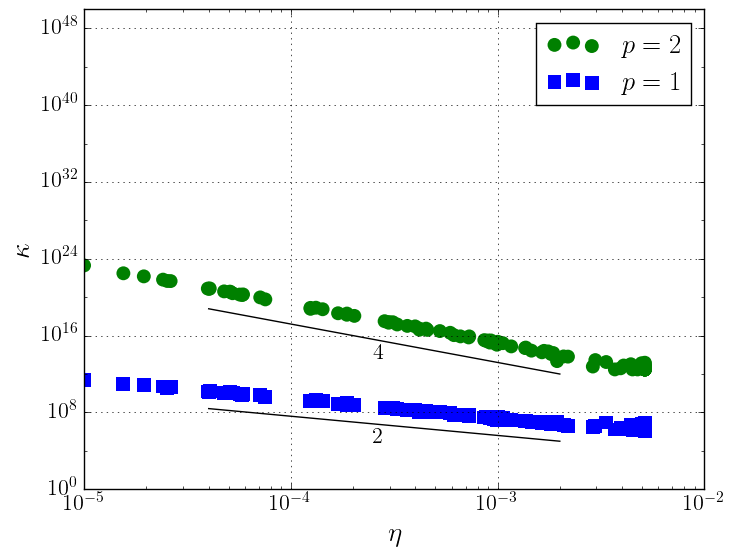}
   \caption{Lagrange}
   \label{fig:condition_lagrange}
 \end{subfigure}
 \caption{Condition number \emph{vs.}\ the smallest volume fraction for the unit square with circular exclusion for (a) B-splines of order $p \in \{ 1,2,3,4\}$, and (b) Lagrange polynomials of order $p \in \{ 1, 2 \}$.}
 \label{fig:condition}
\end{figure}

A power algorithm is used to compute the smallest and largest eigenvalues of the system matrix, and thereby the condition number. The algorithm is terminated when the Rayleigh quotient of two subsequent vectors has a relative difference of less than $10^{-6}$. The inverse power method used to compute the smallest eigenvalue, relies on the inversion of the ill-conditioned finite cell system matrix. The preconditioner developed in this manuscript is used to reliably execute this inverse power iteration. When matrices are singular up to machine precision, the linear dependence tolerance used to construct the preconditioner (see Section~\ref{sec:algorithm}) can hinder the computation of the smallest eigenvalue. As we will study in further detail in the next section, the underlying mechanism leading to ill-conditioning is generally different for B-spline bases and Lagrange bases. As a result of this difference, the condition numbers for B-spline discretizations can be computed accurately over the full range of data shown in Figure~\ref{fig:condition}, while machine precision hinders the computation of condition numbers for Lagrange bases of third and fourth order. Therefore, we do not present results for the Lagrange bases for $p \in \{3,4\}$.

Figure~\ref{fig:condition} shows the results for B-splines of order $p \in \{1,2,3,4\}$ and for Lagrange basis functions of order $p \in \{1,2\}$. The upper bound of the volume fraction $\eta_R=5 \cdot 10^{-3}$ is observed from these plots. The predicted scaling rates of the condition number according to \eqref{eq:estimate} are depicted in black. Evidently, the numerical results for all considered basis functions closely resemble these scaling rates. Also in agreement with the above derivations is the observation that the condition numbers obtained for Lagrange bases are very similar to those obtained for B-splines of equal order. 

\section{Algebraic preconditioning for FCM}
\label{sec:precon}
Motivated by the condition number analysis presented above, in this section we develop an FCM preconditioner. The fundamental observation on which our developments are based is that ill-conditioning of finite cell systems is caused by the occurrence of small eigenvalues of the system matrix $\mathbf{A}$ due to small trimmed cell volume fractions (see \eqref{eq:est_invnorm}). Recall that from equations \eqref{eq:norm1} and \eqref{eq:funcvecratio}, it follows that the smallest eigenvalue for finite cell systems can be expressed as 
\begin{equation}
 \lambda_{\rm min} = \frac{1}{\|\mathbf{A}^{-1}\|_2} =
 \min_{\mathbf{y}} \frac{\|\mathbf{y}\|_{\mathbf{A}}^2}{\|\mathbf{y}\|_2^2} =
 \min_{v_h,\mathbf{y}} \frac{\|v_h\|_{\mathcal{F}}^2}{\| \mathbf{y} \|_2^2},
 \label{eq:lambdamin}
\end{equation}
where the function $v_h \in \mathcal{V}_h$ corresponds to a vector $\mathbf{y} \in \mathbb{R}^n$ by the relation $v_h = \boldsymbol{\Phi}^T\mathbf{y}$,
with $\boldsymbol{\Phi}$ denoting the vector containing all functions in the basis of $\mathcal{V}_h$. Hence, ill-conditioning occurs when the basis allows functions $v_h$ and corresponding vectors $\mathbf{y}$ for which $\|v_h\|_{\mathcal{F}} \ll \|\mathbf{y}\|_2$. To improve the conditioning, we construct an alternative basis, $\overline{\boldsymbol{\Phi}} \, =\mathbf{S}\boldsymbol{\Phi}$, which precludes these large differences between the norm of functions and the norm of the coefficient vectors. For nonsingular preconditioning matrices $\mathbf{S}$, both $\boldsymbol{\Phi}$ and its preconditioned counterpart $\overline{\boldsymbol{\Phi}}$ span the same approximation space $\mathcal{V}_h$.

For every function $v_h \in \mathcal{V}_h$, there are unique coefficient vectors $\mathbf{y}$ and $\overline{\mathbf{y}}$ such that $v_h = \boldsymbol{\Phi}^T\mathbf{y} = \overline{\boldsymbol{\Phi}}^T \, \overline{\mathbf{y}} \, = \boldsymbol{\Phi}^T\mathbf{S}^T\overline{\mathbf{y}}$, from which it follows that $\mathbf{y} = \mathbf{S}^T\overline{\mathbf{y}}$. The preconditioned solution $\overline{\mathbf{x}}$, is obtained by solving the symmetrically preconditioned linear system $\overline{\mathbf{A}} \, \overline{\mathbf{x}} \, = \overline{\mathbf{b}}$, where the preconditioned system matrix $\overline{\mathbf{A}}$ and right hand side vector $\overline{\mathbf{b}}$ are given by
\begin{equation}
\begin{aligned}
 \overline{\mathbf{A}} &= \mathcal{F}(\overline{\boldsymbol{\Phi}},\overline{\boldsymbol{\Phi}}^T) & &= \mathcal{F}(\mathbf{S}\boldsymbol{\Phi},\boldsymbol{\Phi}^T\mathbf{S}^T) & &=  \mathbf{S}\mathcal{F}(\boldsymbol{\Phi},\boldsymbol{\Phi}^T)\mathbf{S}^T & &= \mathbf{S}\mathbf{A}\mathbf{S}^T,\\
 \overline{\mathbf{b}} &= \ell(\overline{\boldsymbol{\Phi}}) & &= \ell(\mathbf{S}\boldsymbol{\Phi}) & &= \mathbf{S}\ell(\boldsymbol{\Phi}) & &= \mathbf{S}\mathbf{b}.
\end{aligned}
\label{eq:sympreconsys}
\end{equation}
In these expressions, $\mathbf{A}$ and $\mathbf{b}$ are the system matrix and right hand side vector corresponding to the original basis $\boldsymbol{\Phi}$. Given an original system $\mathbf{A} \mathbf{x}=\mathbf{b}$ and preconditioning matrix $\mathbf{S}$, the resulting symmetrically preconditioned system is
\begin{equation}\begin{aligned}
\mathbf{S}\mathbf{A}\mathbf{S}^T \overline{\mathbf{x}} & = \mathbf{S}\mathbf{b}, \\
\mathbf{x} & = \mathbf{S}^T \overline{\mathbf{x}}.
\end{aligned}\end{equation}

In the remainder of this section, we will discuss the development of an effective FCM preconditioner. In Section \ref{sec:cond_sources}, we first identify two sources for the occurrence of small eigenvalues. Each source of ill-conditioning can be remedied by a modification of the basis through a preconditioning matrix, which can be constructed algebraically (based solely on the original system matrix). Because this preconditioner is not interwoven with the rest of the method, it is robust and straightforward to automate and implement. Since this preconditioner does not manipulate the weak form problem or approximation space, the obtained solution is unaffected, and stability properties of the employed approximation space (\emph{e.g.,} the inf-sup condition for mixed methods) are maintained. In Section~\ref{sec:sipic}, we will discuss how the developed preconditioner -- which we refer to as \emph{Symmetric Incomplete Permuted Inverse Cholesky} (SIPIC) -- is related to well-established preconditioning strategies, after which we discuss various implementation aspects in Section~\ref{sec:algorithm}. Finally, we demonstrate the effectivity of the SIPIC preconditioner by means of numerical simulations in Section~\ref{sec:preconexample}. 

\subsection{Small eigenvalues: sources and remedies}\label{sec:cond_sources}

\subsubsection{Diagonal scaling}
\label{sec:scaling}
\begin{figure}
  \centering
  \begin{subfigure}[t]{0.49\textwidth}
    \includegraphics[width=\textwidth]{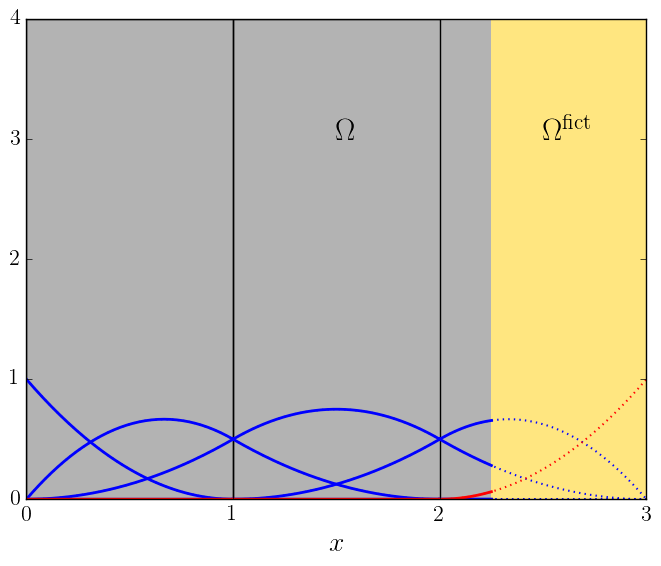}
    \caption{Original}
    \label{fig:spline_orig}
  \end{subfigure}
  \begin{subfigure}[t]{0.49\textwidth}
    \includegraphics[width=\textwidth]{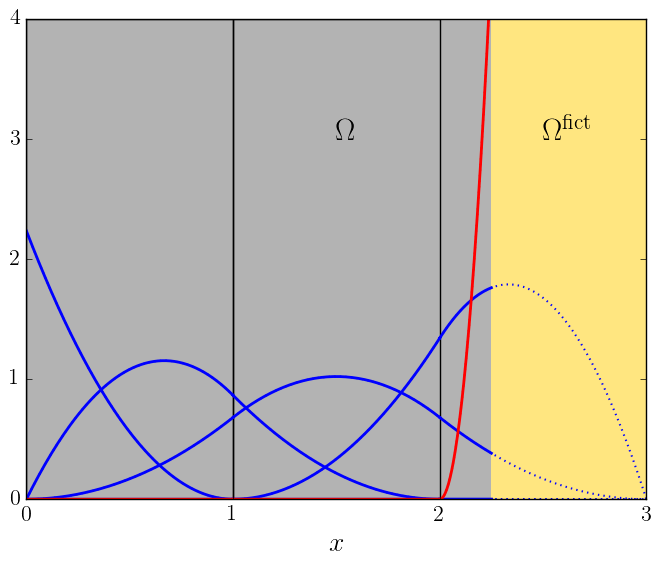}
    \caption{Scaled}
    \label{fig:spline_diag}
  \end{subfigure}
  \caption{Second-order B-spline basis functions over the domain $\Omega=(0,2\frac{1}{4})$ created using the knot vector $[0,0,0,1,2,3,3,3]$: (a) before scaling, and (b) after scaling.}
  \label{fig:scaling_splines}
\end{figure}

A small eigenvalue can be caused by the occurrence of a (trimmed) basis function $\phi$, with a small FCM norm relative to other (possibly untrimmed) basis functions. This can be inferred from equation \eqref{eq:lambdamin} by considering the standard unit vector $\mathbf{y}$ ($\|\mathbf{y}\|_2=1$) corresponding to $v_h=\phi$. This situation is typical for B-spline basis functions whose support contains only a single trimmed cell with a small volume fraction, since all the derivatives up to order $p-1$ of such basis functions vanish on the (original) boundaries of that cell. We illustrate this situation in Figure~\ref{fig:spline_orig}, where the basis function with a small norm (depicted in red) is similar to the function specified in \eqref{eq:smallfunc}.

To remedy this source of small eigenvalues, all basis functions $\phi$ are normalized with respect to the norm $\|\phi\|_{\mathcal{F}}$, which yields the modified basis functions $\phi^* = \phi/\|\phi\|_{\mathcal{F}}$. This diagonal scaling operation can be cast into the form of a preconditioning matrix $\mathbf{D}$ whose diagonal corresponds to the reciprocal of the square root of the main diagonal of system matrix $\mathbf{A}$:
\begin{equation}
\mathbf{D} = \begin{bmatrix} \frac{1}{\sqrt{A_{11}}} & & \\ & \ddots & \\ & & \frac{1}{\sqrt{A_{nn}}} \end{bmatrix} = \begin{bmatrix} \frac{1}{{\|\phi_1\|_{\mathcal{F}}}} & & \\ & \ddots & \\ & & \frac{1}{{\|\phi_n\|_{\mathcal{F}}}} \end{bmatrix}.
\label{eq:scalingmatrix}\end{equation}
In Figure~\ref{fig:spline_diag}, we illustrate the effect of this diagonal scaling operation for the basis $\boldsymbol{\Phi}$ in Figure~\ref{fig:spline_orig}, where the preconditioned basis follows from $\boldsymbol{\Phi}^*= \mathbf{D} \boldsymbol{\Phi}$. We note that diagonal preconditioning has also been demonstrated to be effective in ameliorating similar ill-conditioning problems encountered in XFEM~\cite{DiagonalXFEM}.

\subsubsection{Local orthonormalization}\label{sec:orthonormalization}
Even if basis functions are properly scaled, small eigenvalues can occur when the global linear independence property of two (or more) basis functions is compromised by a trimming operation. This results in the near linear dependence of rows and columns of $\mathbf{A}$ associated with the trimmed cell basis functions and, accordingly, results in ill-conditioning of $\mathbf{A}$. For example, consider two (scaled) basis functions, $\phi_1^*$ and $\phi_2^*$, that are very similar with respect to the finite cell norm, \emph{i.e.,} $\| \phi_1^* - \phi_2^* \|_{\mathcal{F}} \ll 1$ (or $\| \phi_1^* + \phi_2^* \|_{\mathcal{F}} \ll 1$), where use is made of $\| \phi_1^* \|_{\mathcal{F}} = \| \phi_2^* \|_{\mathcal{F}} = 1$. The difference (or sum) of these functions corresponds to a coefficient vector $\mathbf{y}$ with $\|\mathbf{y}\|_2=\sqrt{2}$. From Equation~\eqref{eq:lambdamin} it then follows that the minimal eigenvalue is bounded from above by $\frac{1}{2}\| \phi_1^* - \phi_2^* \|_{\mathcal{F}}^2 \ll 1$ (or $\frac{1}{2}\| \phi_1^* + \phi_2^* \|_{\mathcal{F}}^2 \ll 1$). We refer to this source of ill-conditioning as \emph{quasi linear dependence} of the trimmed basis.

\begin{figure}
  \centering
    \begin{subfigure}[t]{0.49\textwidth}
    \includegraphics[width=\textwidth]{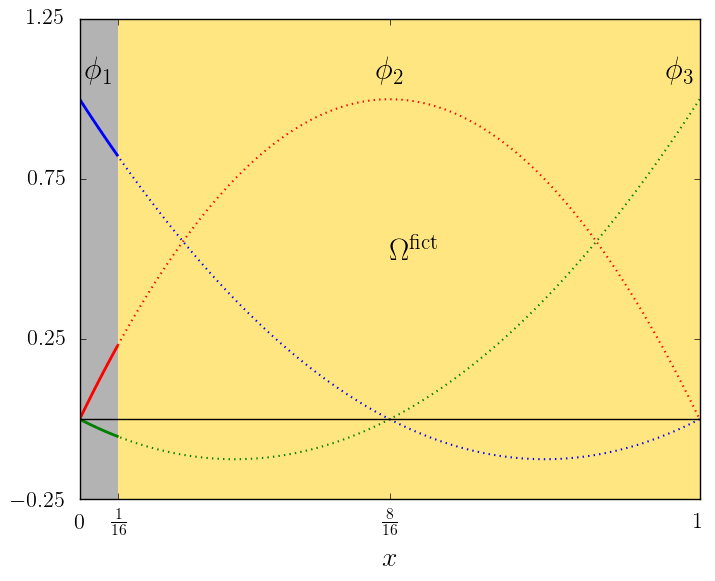}
    \caption{Original on untrimmed cell}
    \label{fig:lagrange_orig}
  \end{subfigure}
  \begin{subfigure}[t]{0.49\textwidth}
    \includegraphics[width=\textwidth]{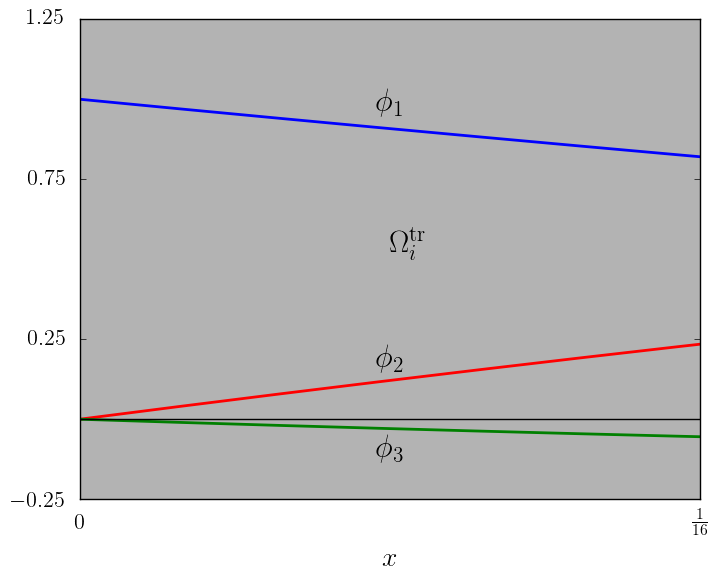}
    \caption{Original on trimmed cell}
    \label{fig:lagrange_restrict}
  \end{subfigure}
  
  \vspace{2em}
  
  \begin{subfigure}[t]{0.49\textwidth}
    \includegraphics[width=\textwidth]{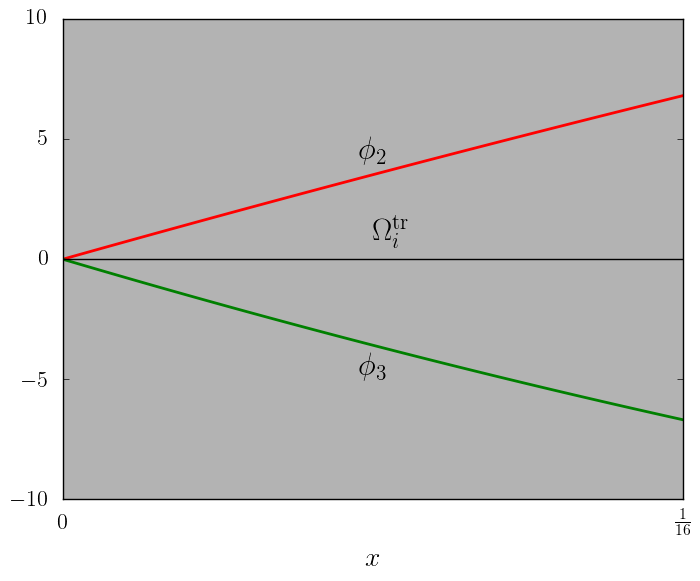}
    \caption{Scaled on trimmed cell}
    \label{fig:lagrange_diag}
  \end{subfigure}  
  \begin{subfigure}[t]{0.49\textwidth}
    \includegraphics[width=\textwidth]{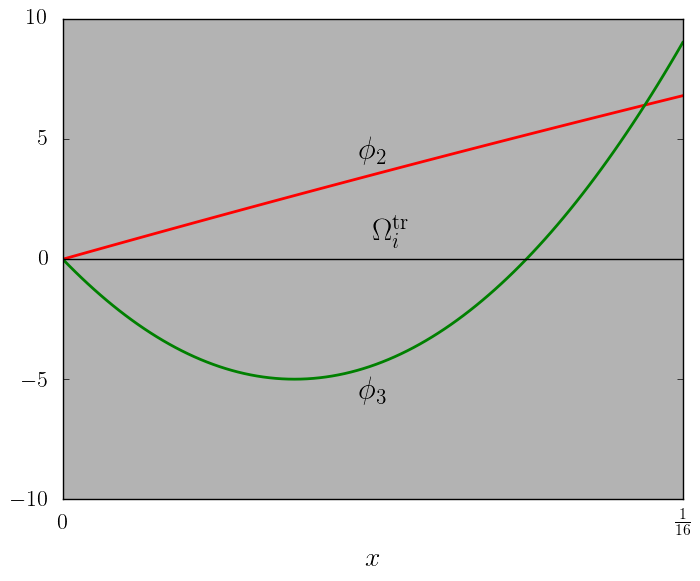}
    \caption{Orthonormalized on trimmed cell}
    \label{fig:lagrange_precon}
  \end{subfigure}
  
  \caption{Second-order Lagrange basis functions over the right-most cell of the discretization with $h=1$ of a trimmed one-dimensional domain. The volume fraction of this cell is equal to $\eta = \frac{1}{16}$.}\label{fig:orthonormalizing_lagrange}
\end{figure}

In Figure~\ref{fig:orthonormalizing_lagrange}, we illustrate the problem of quasi linear dependence for a one dimensional discretization with second order Lagrange basis functions. We consider the rightmost cell, $\Omega_i^{\rm tr} = (\hat{x},\hat{x}+\widetilde{h})=(0,\frac{1}{16})$, of a discretization with mesh size $h=1$ which is trimmed by a natural boundary (\emph{i.e.,} no boundary terms). While the support of basis function $\phi_1$ extends beyond the trimmed element, the basis functions $\phi_2$ and $\phi_3$, 
\begin{equation}\begin{aligned}
\phi_2  &= 4 \left(\frac{x-\hat{x}}{h}\right) - 4 \left(\frac{x-\hat{x}}{h}\right)^2 ,\\
\phi_3   &=   - \left(\frac{x-\hat{x}}{h}\right) + 2 \left(\frac{x-\hat{x}}{h}\right)^2,
\end{aligned}\end{equation}
are only supported on this cell. On the untrimmed cell, these basis functions are evidently linearly independent. For instance, the $\mathcal{F}$-angle between $\phi_2$ and $\phi_3$ is
\begin{equation}
 |\cos(\theta_{\mathcal{F}})| = \frac{|\mathcal{F}(\phi_2,\phi_3)|}{\| \phi_2 \|_{\mathcal{F}} \| \phi_3 \|_{\mathcal{F}}} = \sqrt{4/7},
\end{equation}
which is (in absolute sense) significantly smaller than unity. However, when $\eta=\widetilde{h}/h \ll 1$, the the $\mathcal{F}$-angle is
\begin{equation}
 |\cos(\theta_{\mathcal{F}})| = \frac{|-1 + 3 \eta - \frac{8}{3} \eta^2|}{\sqrt{1-6\eta+\frac{44}{3}\eta^2-16\eta^3+\frac{64}{9}\eta^4}} \approx 1 - \frac{1}{6}\eta^2,
\end{equation}
which reveals a quasi linear dependence. Similar observations can be made when considering linear combinations of these functions. In the limit of $\eta$ going to zero, the local coordinate $(x-\hat{x})/h \ll 1$, which reduces both $\phi_2$ and $\phi_3$ to linears on $\Omega_i^{\rm tr}$. Figure~\ref{fig:lagrange_diag} shows the scaled basis functions for $\eta=\frac{1}{16}$,
\begin{equation}\begin{aligned}
{\phi}_2^* & = \frac{\phi_2}{\|\phi_2\|_{\mathcal{F}}} \approx \sqrt{\frac{h}{\eta}} \left(\frac{x-\hat{x}}{h}\right) -  \sqrt{\frac{h}{\eta}} \left(\frac{x-\hat{x}}{h}\right)^2,\\
{\phi}_3^*  & = \frac{\phi_3 }{\|\phi_3 \|_{\mathcal{F}}} \approx -  \sqrt{\frac{h}{\eta}}\left(\frac{x-\hat{x}}{h}\right) + 2  \sqrt{\frac{h}{\eta}} \left(\frac{x-\hat{x}}{h}\right)^2,
\end{aligned}\end{equation}
where it is used that $\|\phi_2 \|_{\mathcal{F}} \approx 4 \sqrt{\eta/h}$ and $\|\phi_3 \|_{\mathcal{F}} \approx \sqrt{\eta/h}$ if $\eta \ll 1$. When these functions are added, the linear term cancels and a function of the form \eqref{eq:smallfunc} results:
\begin{equation}
{\phi}_2^* + {\phi}_3^* \approx \sqrt{\frac{h}{\eta}} \left(\frac{x-\hat{x}}{h}\right)^2,
\end{equation}
for which it can be shown that $\|{\phi}_2^* + {\phi}_3^*\|_{\mathcal{F}} \approx \frac{4}{3} \eta^2 \ll 1$. Hence, the scaled basis functions $\phi_2^*$ and $\phi_3^*$ are quasi linearly dependent for small volume fractions $\eta$. In this case, diagonal scaling is not effective in the sense that the upper bound of the minimal eigenvalue \eqref{eq:lambdamin}, and thereby the condition number, remains dependent on the volume fraction $\eta$, leading to ill-conditioning for small volume fractions. It is noted, however, that the scaling rate is improved with respect to that of the original system (see \eqref{eq:estimate} and Section~\ref{sec:preconexample}).

To remedy ill-conditioning due to quasi linear dependence, we propose to orthonormalize the quasi linearly dependent functions by the Gram-Schmidt procedure (see \cite{Golub2012} for details). To illustrate how this Gram-Schmidt procedure can be used to construct a preconditioner, we again consider two quasi linearly dependent scaled functions, $\phi_\alpha^*$ and $\phi_\beta^*$. To improve conditioning, the function $\phi_\beta^*$ is made orthogonal to $\phi_\alpha^*$ with respect to the finite cell inner product:
\begin{equation}
 \phi_\beta^{\perp} = \phi_\beta^* - \frac{(\phi_\beta^*,\phi_\alpha^*)_{\mathcal{F}}}{\| \phi_\alpha^* \|_{\mathcal{F}}^2} \phi_\alpha^* = \phi_\beta^* - (\phi_\beta^*,\phi_\alpha^*)_{\mathcal{F}} \phi_\alpha^*.
\end{equation}
In terms of the scaled basis $\boldsymbol{\Phi}^*$, the orthonormalization procedure can be expressed as a matrix-vector operation:
\begin{equation}
 \boldsymbol{\Phi}^\perp = \mathbf{G} \boldsymbol{\Phi}^* = \mathbf{G} \mathbf{D} \boldsymbol{\Phi},
 \label{eq:orthooperation}
\end{equation}
where $\mathbf{G}$ contains zeros everywhere, except for the entries
\begin{equation}
   G_{\beta\alpha} =    - (\phi_\beta^*,\phi_\alpha^*)_{\mathcal{F}} = - [\mathbf{DAD}^T]_{\beta\alpha} = - \frac{A_{\beta\alpha}}{\sqrt{A_{\beta\beta}A_{\alpha\alpha}}}      \quad \mbox{and} \quad     G_{\gamma\gamma} = 1.
\end{equation}
Since the norm of the basis functions is not preserved by the orthogonalization operation \eqref{eq:orthooperation}, in order to avoid improperly scaled basis functions leading to ill-conditioning, the orthogonalization procedure is followed by another normalization operation:
\begin{equation}
 \overline{\boldsymbol{\Phi}} = \mathbf{D}^\perp \mathbf{G} \boldsymbol{\Phi}^* = \mathbf{D}^\perp \mathbf{G} \mathbf{D} \boldsymbol{\Phi},
\end{equation}
where $\mathbf{D}^\perp$ is the scaling matrix corresponding to $\mathbf{A}^\perp= \mathbf{G} \mathbf{A}^* \mathbf{G}^T=\mathbf{G} \mathbf{D} \mathbf{A} \mathbf{D}^T\mathbf{G}^T$. The result of this orthonormalization procedure is the preconditioning matrix $\mathbf{S}=\mathbf{D}^\perp \mathbf{G} \mathbf{D}$, which is algebraic by virtue of the fact that the three matrices $\mathbf{D}^\perp$, $\mathbf{G}$ and $\mathbf{D}$ can all be constructed on the basis of the original system matrix $\mathbf{A}$. We note that similar procedures have been applied in the context of SGFEM in order to orthonormalize the additional functions to the existing FEM approximation space~\cite{Babuska2012}.

With higher order methods or in more than one dimension, it can occur that more than two functions are quasi linearly dependent, which requires multiple orthonormalization steps in the Gram-Schmidt procedure. The above concept can be applied without fundamental modifications in such cases, which is discussed in Section~\ref{sec:algorithm}. In the absence of boundary contributions to the system matrix, quasi linear dependence predominantly occurs between functions with the same support in the physical domain. The orthonormalization procedure outlined above does not change the support of basis functions in that case, and preserves the sparsity pattern of the system matrix. In the case of Dirichlet boundaries, it can incidentally occur that a cell is trimmed in such a way that  the penalty term dominates the bilinear form for all functions supported on that cell, including functions that are also supported on other cells in the physical domain. In these cases, the Gram-Schmidt procedure may increase the support of basis functions, causing additional fill-in in the preconditioned system matrix. In our algorithm we limit this fill-in by appropriately ordering the quasi linearly dependent basis functions (see Section~\ref{sec:algorithm}). The effect of fill-in is studied numerically in Section~\ref{sec:preconexample}.

\subsection{Preconditioner characterization and relation to alternative techniques}
\label{sec:sipic}
The preconditioner $\mathbf{S}$ developed herein approximates the inverse of the system matrix by:
\begin{equation}
 \mathbf{S} \mathbf{A} \mathbf{S}^T \approx \mathbf{I} \qquad\mbox{or}\qquad \mathbf{S}^{T}  \mathbf{S} \approx \mathbf{A}^{-1}.
 \label{eq:inverseapprox}
\end{equation}
These identities would be exact if the Gram-Schmidt procedure would be applied to all basis functions. The developed preconditioner is incomplete in the sense that only quasi linearly dependent basis functions are orthonormalized. By construction, there exists a permutation matrix $\mathbf{P}$ which reorders the rows of $\mathbf{S}$ in such a way that a lower triangular matrix is obtained: $\mathbf{L}=\mathbf{P}\mathbf{S}\mathbf{P}^T$. Substitution in \eqref{eq:inverseapprox} yields
\begin{equation}
\mathbf{L}^T \mathbf{L} \approx \mathbf{P}  \mathbf{A}^{-1}  \mathbf{P}^T = [ \mathbf{P}\mathbf{A}\mathbf{P}^T ]^{-1} =\mathbf{C}^{-T} \mathbf{C}^{-1},
\label{eq:LTL}
\end{equation}
where use has been made of the property $\mathbf{P}^{T}  \mathbf{P}=\mathbf{I}$ of permutation matrices and $\mathbf{C}$ corresponds to the Cholesky decomposition of the permuted matrix, \emph{i.e.,} $\mathbf{C} \mathbf{C}^{T} = \mathbf{P}\mathbf{A}\mathbf{P}^T $. Equation \eqref{eq:LTL} conveys that the lower triangular matrix $\mathbf{L}$ approximates the inverse Cholesky decomposition of the permuted matrix, $\mathbf{L} \approx \mathbf{C}^{-1}$, and that the preconditioner can be written as $\mathbf{S} = \mathbf{P}^T \mathbf{L} \mathbf{P} \approx \mathbf{P}^T \mathbf{C}^{-1} \mathbf{P}$. Hence this preconditioner can be interpreted as the conjugate permutation of an incomplete inverse Cholesky decomposition of a permutation of the system matrix. Note that this permutation is done to reduce fill-in as described in the last paragraph of Section~\ref{sec:orthonormalization}. Since we apply $\mathbf{S}$ symmetrically (see Equation~\eqref{eq:sympreconsys}) we refer to this technique as \emph{Symmetric Incomplete Permuted Inverse Cholesky} (SIPIC) preconditioning.

\subsection{Preconditioner-construction algorithm}
\label{sec:algorithm}
Algorithm~\ref{alg:precon} outlines the construction of the SIPIC preconditioner $\mathbf{S}$. Since this preconditioner is algebraic, it is constructed solely based on the information contained in the system matrix $\mathbf{A}$. As the SIPIC preconditioner is incomplete in the sense that only functions that are quasi linearly dependent are orthonormalized, an orthonormalization threshold parameter $\gamma\in[0,1]$ is to be provided as input to the algorithm.

The algorithm is initialized with the construction of the scaling matrix $\mathbf{D}$ and the sets $\mathcal{I}$ and $\mathcal{J}$ that respectively identify the current and total quasi linear dependencies. Subsequently, the linear dependencies are grouped in $\Sigma$, and for every sorted local group of indices $\sigma \in \Sigma$ the Gram-Schmidt orthonormalization procedure is applied. Finally, a check is performed to test if any new quasi linear dependencies have emerged after this orthonormalization procedure. If this is the case, these are added to the total set of quasi linear dependencies and the orthonormalization loop is restarted. Otherwise the construction of the SIPIC preconditioner is finished. In our simulations, we have observed convergence of this orthonormalization procedure in at most two iterations, resulting in negligible computational cost for the construction of the SIPIC preconditioner.

The subroutines in Algorithm~\ref{alg:precon} perform the following operations:
\begin{itemize}
\item \verb;scale(;$\mathbf{A}$\verb;); returns the diagonal matrix $\mathbf{D}$ according to Equation~\eqref{eq:scalingmatrix}.
\item \verb;identify(;$\mathbf{SAS}^T,\gamma$\verb;); locates quasi linear dependencies in the scaled system $\mathbf{SAS}^T$ (\emph{i.e.,} the diagonal entries are equal to one). A measure for the linear dependence of two basis functions with indices $\alpha$ and $\beta$ is derived from the Cauchy-Schwarz inequality:
\begin{equation}
 \frac{|(\phi_\alpha,\phi_\beta)_{\mathcal{F}}|}{\|\phi_\alpha\|_{\mathcal{F}}\|\phi_\beta\|_{\mathcal{F}}}  \leq 1,
\end{equation}
where equality to one indicates linear dependence. Two basis functions are identified as being quasi linearly dependent if the absolute value of the corresponding off-diagonal term in the scaled system matrix exceeds the orthogonalization threshold, 
\begin{equation}
 \left| [\mathbf{SAS}^T]_{\alpha\beta}  \right| > \gamma,
 \label{eq:diagthres}
\end{equation}
which -- exploiting the symmetry of the system matrix -- results in the set of index pairs of quasi linearly dependent basis functions:
\begin{equation}
 \mathcal{I} = \left\{ (\alpha,\beta) \mid \alpha > \beta,~\left| [\mathbf{SAS}^T]_{\alpha\beta} \right|  > \gamma \right\}.
\end{equation}
For all our simulations, a value of $\gamma = 0.9$ was found to be adequate. Larger values may omit some quasi linear dependencies, resulting in a less effective preconditioner. Smaller values can tag functions with different supports as being quasi linearly dependent, causing additional fill-in in the preconditioned system. We study this fill-in behavior in more detail in the context of numerical simulations in Section~\ref{sec:preconexample}. The observed robustness and effectivity of $\gamma=0.9$ does not necessarily extend to other (non-polynomial) bases or applications (such as XFEM or other immersed techniques) suffering from similar conditioning problems. Hence, careful selection of $\gamma$ is required when applying the SIPIC preconditioner in these situations.

\item \verb;group(;$\mathcal{I},\mathbf{A}$\verb;); returns a set ($\Sigma$) of non-intersecting sets ($\sigma\in\Sigma$) in which intersecting sets in $\mathcal{I}$ are replaced by their union. That is, an element $\sigma \in \Sigma$ is a tuple containing the indices of all basis functions that are quasi linear dependent to each other. This tuple is ordered from the index whose corresponding row in the system matrix has the least nonzero entries to the index with the largest number of nonzero row entries in $\mathbf{A}$, in order to reduce fill-in as described in the last paragraph of Section~\ref{sec:orthonormalization}.
\item \verb;orthonormalize(;$\mathbf{A}_\sigma$\verb;); applies the Gram-Schmidt orthonormalization to the functions indexed by $\sigma$, such that $\mathbf{S}_\sigma$ becomes the inverse Cholesky decomposition of $\mathbf{A}_\sigma$ and $\mathbf{S}_\sigma\mathbf{A}_\sigma\mathbf{S}_\sigma^T$ becomes the identity matrix.

Instabilities of the construction algorithm can occur when -- due to the presence of cells with extremely small volume fractions -- some basis functions become linearly dependent in the sense that the absolute off-diagonal entries \eqref{eq:diagthres} are up to machine precision equal to one. When this occurs, $[\mathbf{SAS}^T]_{ii} = 0$ for some index $i$ after orthogonalization, and rescaling will result in a division by zero. To stabilize the algorithm, a check is performed to test if $[\mathbf{SAS}^T]_{ii} > \varepsilon$ before rescaling. If this is not the case, the function is eliminated by simply deleting the row $i$ from the system. The preconditioner $\mathbf{S}$ is then no longer square, but $\mathbf{SAS}^T$ remains square and SPD and is merely reduced in size. This operation does not influence the quality of the obtained solution, because the remaining basis functions span the same approximation space up to machine precision. A value of $\varepsilon = 10^2 \cdot {\tt eps}$ (with ${\tt eps}$ the machine precision) was found to be adequate for all simulations considered in this manuscript. An example of a matrix that is singular up to machine precision as a result of quasi linear dependence is given in \ref{sec:app_singular}.
\end{itemize}

\begin{algorithm}
 \DontPrintSemicolon
 \SetKwFunction{identify}{identify}
 \SetKwFunction{group}{group}
 \SetKwFunction{scale}{scale}
 \SetKwFunction{orthonormalize}{orthonormalize}
 \SetKwFunction{Range}{range}
 \SetKwFunction{Len}{len}
 \SetKwFunction{Call}{}
 
 \SetKwComment{Comment}{\#\,}{}
 \SetKwComment{tcc}{\#\,}{}
 
 \SetStartEndCondition{ }{}{}
 \SetKwProg{Fn}{def}{\string:}{end}
 \SetKw{KwTo}{in}
 \SetKwFor{For}{for}{\string:}{end}
 \SetKwFor{While}{while}{:}{end}

 \KwIn{ $\mathbf{A}$, $\gamma$ {\Comment*[r]{system matrix, orthonormalization threshold}}}
 \KwOut{$\mathbf{S}$ {\Comment*[r]{SIPIC preconditioner}}}
 
 \BlankLine
 
 \Comment*[l]{Initialize}
 $\mathbf{S} =$ \scale{$\mathbf{A}$} \Comment*[r]{diagonal scaling}
 $\mathcal{I} = \mathcal{J} =$ \identify{$\mathbf{SAS}^T,\gamma$} \Comment*[r]{initial dependencies}
 
 \BlankLine
 
 \Comment{Main loop}
 \While{$\mathcal{J} \neq \emptyset$}
 {
  $\Sigma =$ \group{$\mathcal{I},\mathbf{A}$}  \Comment*[r]{group and order dependencies}
  \For{$\sigma$ \KwTo $\Sigma$}{
   $\mathbf{S}$\FuncSty{[}$\sigma,\sigma$\FuncSty{]}$ = $\orthonormalize{$\mathbf{A}$\FuncSty{[}$\sigma,\sigma$\FuncSty{]}}\;
                         }
  $\mathcal{J} =$ \identify{$\mathbf{SAS}^T,\gamma$}\;
  $\mathcal{I} = \mathcal{I} \cup \mathcal{J}$ \Comment*[r]{append current dependencies}
 }
 
 \BlankLine
 
 \Return{$\mathbf{S}$}
 
 \BlankLine
 \BlankLine
 \BlankLine
 \BlankLine

 \Comment*[l]{Orthonormalization function}
 \Fn{\orthonormalize{$\mathbf{A}_\sigma$}}
 {
 
  \BlankLine
  
  $\mathbf{S}_\sigma =$ \scale{$\mathbf{A}_\sigma$} \Comment*[r]{diagonal scaling}
  
  \BlankLine
  
  \For(\tcc*[f]{Gram-Schmidt orthonormalization}){$i$ \KwTo \Range{\Len{$\mathbf{S}_\sigma$}}}{
   \For(\tcc*[f]{$j<i$}){$j$ \KwTo \Range{$i$}}{
    $\mathbf{S}_\sigma\FuncSty{[}i,:\FuncSty{]} = \mathbf{S}_\sigma$\FuncSty{[}$i,:$\FuncSty{]}$ - \FuncSty{(} \mathbf{S}_\sigma\mathbf{A}_\sigma\mathbf{S}_\sigma^T\FuncSty{)}\FuncSty{[}i,j\FuncSty{]} \mathbf{S}_\sigma \FuncSty{[}j,:\FuncSty{]}$ 
              }
   $\mathbf{S}_\sigma\FuncSty{[}i,:\FuncSty{]} = \frac{\mathbf{S}_\sigma \FuncSty{[} i,:\FuncSty{]}}{\FuncSty{(}\sqrt{\mathbf{S}_\sigma\mathbf{A}_\sigma\mathbf{S}_\sigma^T} \FuncSty{)}\FuncSty{[}i,i\FuncSty{]}}$ {\Comment*[r]{rescale}}
  }
           
  \BlankLine
 
  \Return{$\mathbf{S_\sigma}$}
  
  \BlankLine
  
  }
  
 \caption{Construction of the SIPIC preconditioner.}
 \label{alg:precon}
\end{algorithm}

\subsection{Numerical test case continued: Unit square with circular exclusion}
\label{sec:preconexample}
To test the effectivity of the SIPIC preconditioner, we again consider the rotating unit square with circular exclusion as introduced in Section~\ref{sec:cond_example}. Besides the original and the SIPIC-preconditioned condition numbers, also the condition number after diagonal scaling is computed for second order ($p=2$) B-spline and Lagrange bases. For the SIPIC preconditioning, an orthonormalization threshold of $\gamma = 0.9$ is applied.

\begin{figure}
 \centering
 \begin{subfigure}{0.49\linewidth}
   \centering
   \includegraphics[width=\linewidth]{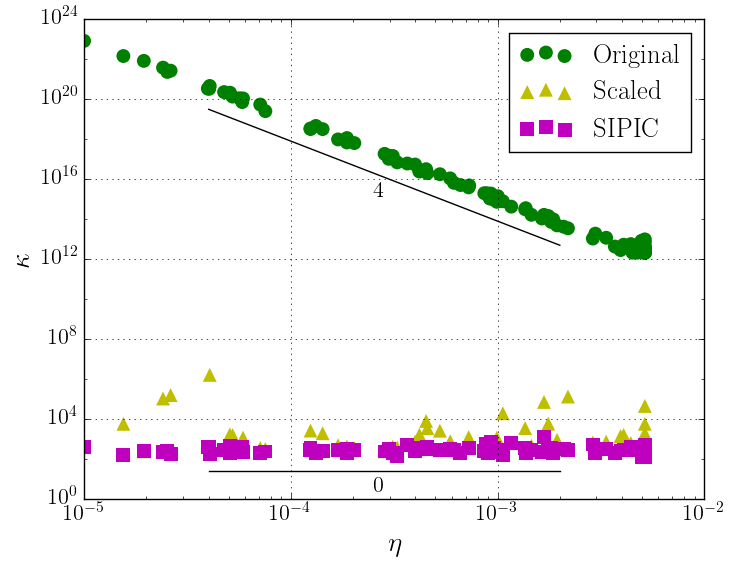}
   \caption{B-splines}
   \label{fig:preconditioned_spline}
 \end{subfigure}
 \begin{subfigure}{0.49\linewidth}
   \centering
   \includegraphics[width=\linewidth]{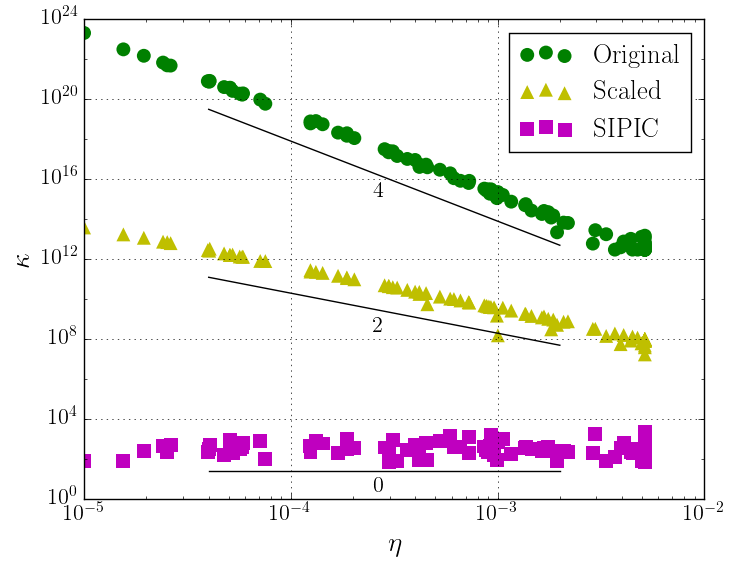}
   \caption{Lagrange}
   \label{fig:preconditioned_lagrange}
 \end{subfigure}
 \caption{Original, scaled and SIPIC-preconditioned condition number \emph{vs.}\ the smallest volume fraction for the unit square with circular exclusion with second order bases.}
 \label{fig:preconditioned}
\end{figure}

The results of all simulations are presented in Figure~\ref{fig:preconditioned}. It is observed that the SIPIC preconditioner drastically improves the condition number for all cases and yields condition numbers independent of the volume fraction (in the sense that there is no scaling relation). In the case of B-splines, diagonal scaling is observed to be effective for most configurations. This is explained by the fact that quasi linear dependencies in B-spline bases are uncommon, since on most of the trimmed cells there is only a single B-spline basis function whose support is restricted to that trimmed cell (see Section~\ref{sec:scaling}). This is in contrast to Lagrange bases, for which quasi linear dependencies are combinatorially more frequent (see Section~\ref{sec:orthonormalization}). As derived for the one-dimensional example in Section~\ref{sec:orthonormalization}, the condition number of the Lagrange discretizations after scaling still scales with $\eta$, but with a slope that is gentler than that of the original condition number. The difference in slope is explained by the fact that all Lagrange basis functions have non-zero first order derivatives on the boundary of their untrimmed support. The basis functions that are only supported on the cell with the smallest volume fraction therefore have a norm of order $\|\phi\|_{\mathcal{F}}^2 \sim \eta^{3-2/d} = \eta^2$. Diagonal scaling consequently increases the magnitude of these functions and the function corresponding to the smallest eigenvalue by a factor $\eta^{-1}$, which explains the observed slope difference of $\eta^{2}$. When a global stabilization parameter is applied, the preconditioner also drastically reduces the condition number. It is observed that -- although not presented here -- the preconditioned condition number still scales with $\eta^{-1/d}$ however, which is exactly the difference in slope between the original condition numbers of locally (see \eqref{eq:estimate}) and globally (see \eqref{eq:estimate_global}) stabilized systems. The origin of this effect and methods to resolve this are beyond the scope of this work. 

\begin{figure}
 \centering
 \begin{subfigure}{0.49\linewidth}
   \centering
   \includegraphics[width=\linewidth]{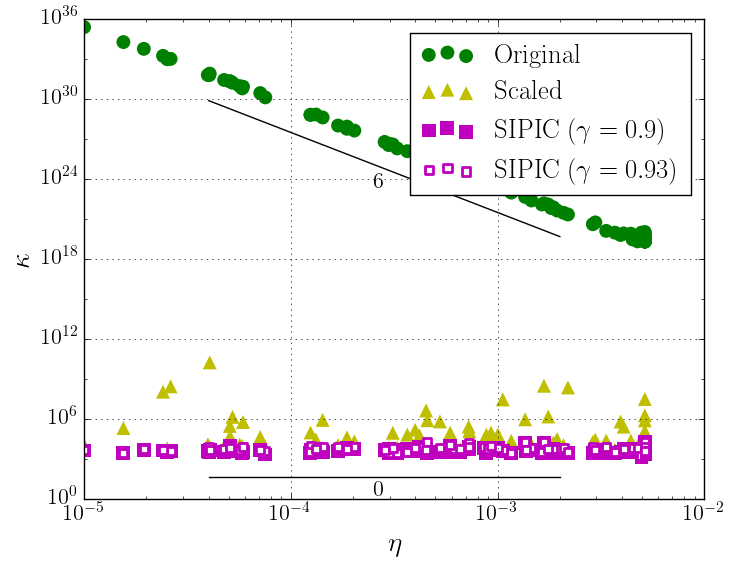}
   \caption{Condition number}
   \label{fig:preconditioned_spline_3}
 \end{subfigure}
 \begin{subfigure}{0.49\linewidth}
   \centering
   \includegraphics[width=\linewidth]{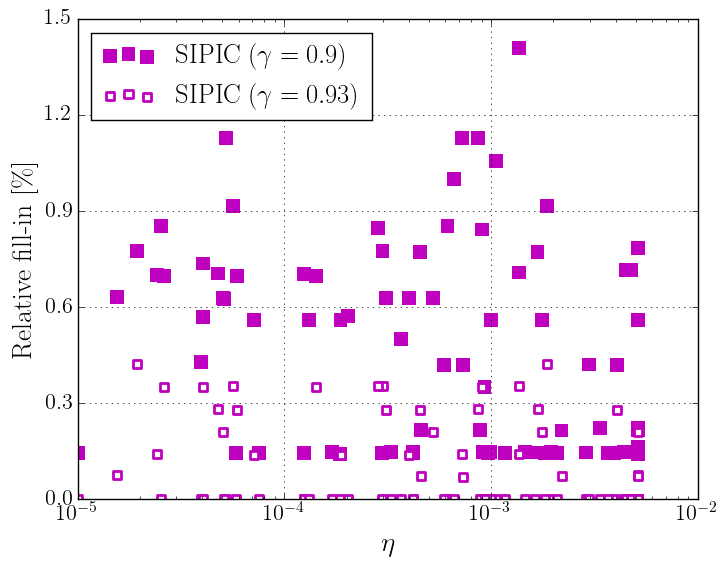}
   \caption{Relative fill-in}
   \label{fig:fillin_spline_3}
 \end{subfigure}
 \caption{Condition number and relative matrix fill-in \emph{vs.}\ the smallest volume fraction for the unit square with circular exclusion with cubic B-splines.}
 \label{fig:spline_3}
\end{figure}

Figure~\ref{fig:fillin_spline_3} presents the system matrix fill-in caused by large boundary contributions (see Section~\ref{sec:orthonormalization}) for two orthonormalization thresholds: $\gamma_1 = 0.9$ and $\gamma_2 = 0.93$. Since no fill-in was observed for second order B-spline bases, we present results for cubic B-splines. The corresponding dependence of the condition number on the smallest volume fraction is shown in Figure~\ref{fig:preconditioned_spline_3}. The relative fill-in is defined as the number of additional nonzero entries of the system matrix divided by the original number of nonzero entries. For $\gamma_1=0.9$, it is observed that fill-in remains limited to approximately $1.5\%$, and that it does not scale with the smallest volume fraction $\eta$. The fill-in can be reduced further at the expense of performing fewer orthonormalization operations by increasing the orthogonalization threshold (see Section~\ref{sec:algorithm}). This is indeed observed from the results for $\gamma_2=0.93$, where fill-in remains limited to approximately $0.5\%$. Although fewer orthonormalization operations are performed, it is observed that the effectivity of SIPIC preconditioning remains unaffected. 

We note that the necessity to increase the threshold is debatable for two reasons:
\begin{itemize}
 \item The fill-in is small and locally contained. An investigation of the functions whose support increases because of the preconditioner conveys that these are all functions for which the bilinear operator is governed by the penalty term on a small cell. Therefore, the effect is confined to functions that are originally supported on that small cell, and cannot spread through the system.
 \item The fill-in is expected to be a two-dimensional artifact. Under the assumption of mesh regularity, it can be shown that the penalty parameter $\beta_i$ is of the order $\beta_i \sim \eta_i^{-1/d}h^{-1}$. Subject to the same assumption, the measure of the edges in a two-dimensional problem is of order $|\Gamma_i| \sim \eta_i^{1/d}h$, such that the total penalty term is of order $\beta |\Gamma_i| \sim 1$. Therefore a small cell (with a limited volumetric contribution) can still have a large penalty term. In the three-dimensional case, the size of the edge is of order $|\Gamma_i| \sim \eta_i^{2/d}h^2$, such that the total penalty term is of order $\beta |\Gamma_i| \sim \eta_i^{1/d}h$. Under the assumption of shape regularity, small cells will therefore have small penalty terms as well, if the number of dimensions is larger than two. As a result, it is unlikely that the penalty term on a cell with a small volume fraction dominates the bilinear form for functions that are also supported on other cells in more than two dimensions.
\end{itemize}
A detailed study of the sensitivity of the SIPIC preconditioning technique to the orthonormalization threshold is a topic of further study.

\section{Numerical simulations: SIPIC preconditioned finite cell analysis }
\label{sec:numerical}
The numerical simulations presented in the previous sections have demonstrated the effectivity of the SIPIC preconditioning technique for the model problem introduced in Section~\ref{sec:FCM}. While these simulations had a strong focus on studying the condition number, in this section we will consider two numerical test cases to demonstrate that SIPIC preconditioning enables robust finite cell analyses using iterative solvers.

Since the considered finite cell systems are symmetric positive definite, a Conjugate Gradient (CG) solver is used. As discussed in Section~\ref{sec:conditioning}, preconditioning of the finite cell system results in faster convergence of the solver (see \eqref{eq:cgconvergence}). Moreover, for well-conditioned systems the residual provides an adequate measure for the error in the energy norm (see \eqref{eq:cgerrorbounds}). As a result, for the same CG tolerance, the energy norm error of a well-conditioned system is generally expected to be smaller than that of an ill-conditioned system. These two beneficial effects of preconditioning enable the computation of high-accuracy finite cell approximations.

For the simulations presented in this section, the solver performance and accuracy aspects of finite cell preconditioning are discussed in the context of mesh convergence studies. Section~\ref{sec:cutplate} considers the elastic analysis of a plate with a circular hole. Since condition numbers can be computed for this two-dimensional test case, it is also used to demonstrate that the results presented in the previous sections extend to linear elasticity. Section~\ref{sec:quarter_ring} studies the performance of SIPIC preconditioning for a three-dimensional linear elasticity problem.

\subsection{Plate with a circular hole}
\label{sec:cutplate}
\begin{figure}
 \centering
 \begin{subfigure}[t]{.49\textwidth}
  \includegraphics[width=\textwidth]{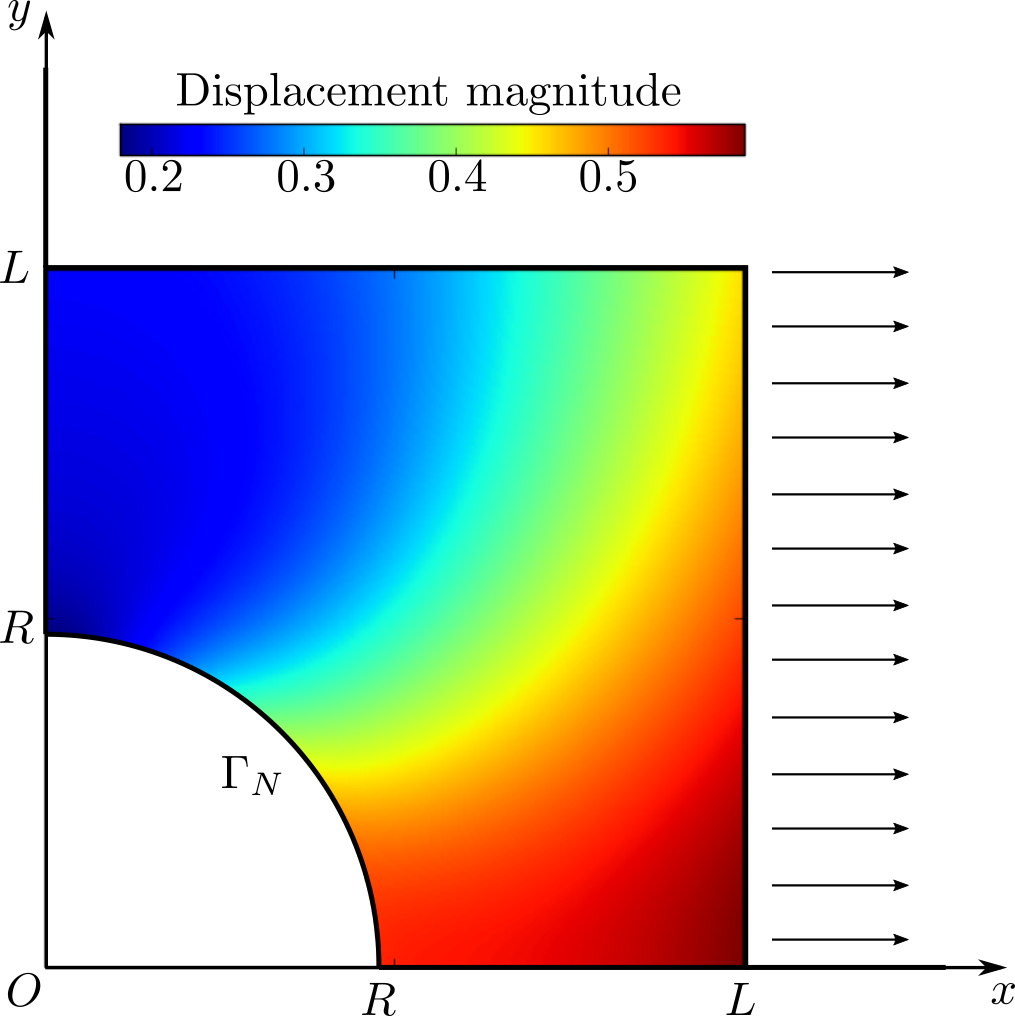}
  \caption{Analytical solution}
  \label{fig:cutplate_solution}
 \end{subfigure}
 \hfill
 \begin{subfigure}[t]{.4\textwidth}
  \includegraphics[width=\textwidth]{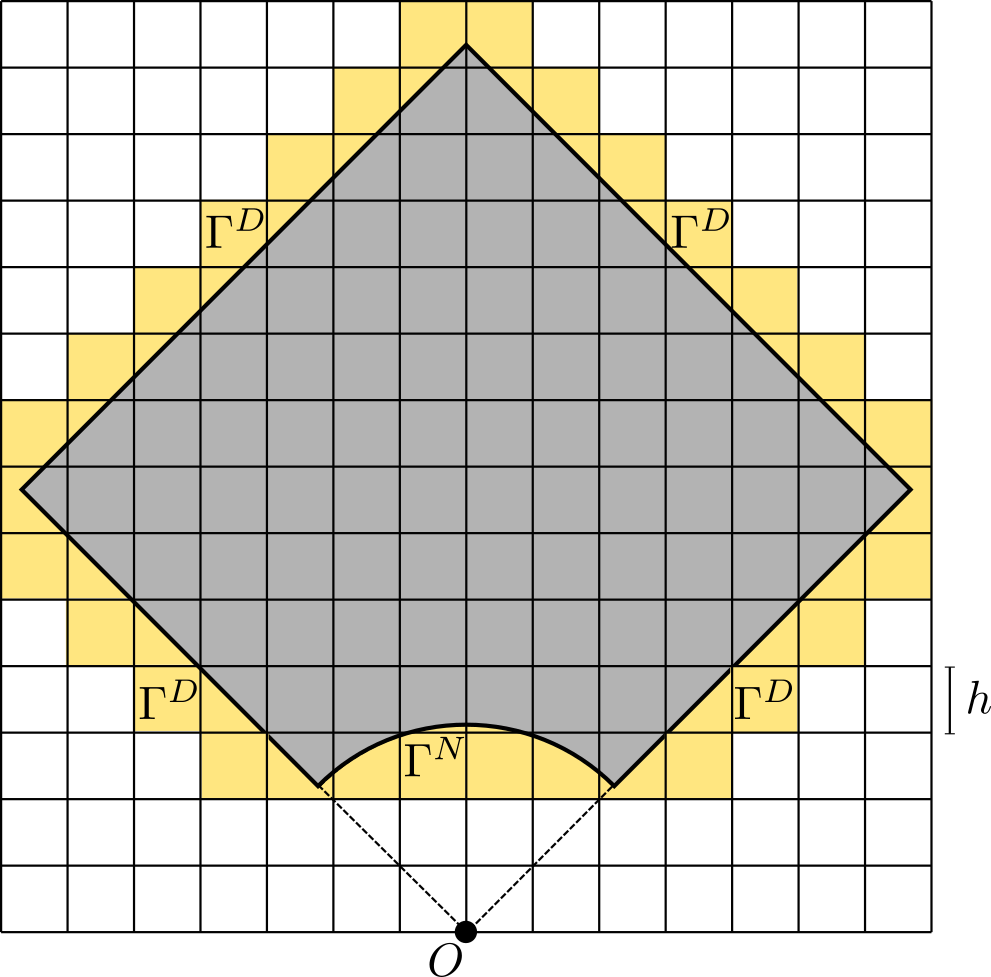}
  \caption{Mesh orientation}
  \label{fig:cutplate_grid}
 \end{subfigure}
 \caption{A plate with a circular hole is rotated over $\SI{45}{\degree}$ with respect to a Cartesian mesh. The Neumann and Dirichlet data on $\Gamma^N$ and $\Gamma^D$, respectively, match the exact solution of the infinite plate problem.}\label{fig:cutplate_testcase}
\end{figure}

We consider the linear elastic analysis of a uniaxially-loaded infinite plate with a circular hole with radius $R=3/(2\pi)$ (Figure~\ref{fig:cutplate_solution}). Under plane strain conditions, the displacement field $u=(u_x,u_y)$ corresponding to a unit horizontal traction loading is given by
\begin{subequations}\begin{align}
u_x & = \frac{x}{\mu} \left( \frac{(2\mu+\lambda)r^2 + (\mu -  \lambda)R^2}{4(\mu+\lambda)r^2} + \frac{\frac{3}{4}R^4 + x^2R^2}{r^4} - \frac{x^2R^4}{r^6} \right),\\
u_y & = \frac{y}{\mu} \left( \frac{    -\lambda  r^2 + (\mu + 3\lambda)R^2}{4(\mu+\lambda)r^2} - \frac{\frac{3}{4}R^4 + y^2R^2}{r^4} + \frac{y^2R^4}{r^6} \right),
\end{align}\label{eq:an_cutplate}\end{subequations}
with $r = \sqrt{x^2 + y^2}$ and Lam\'{e} parameters $\lambda = \mu = 1$. This exact solution is reproduced on a truncated domain by considering the strong formulation
\begin{equation}
  \begin{cases}
    -{\rm div} (\boldsymbol{\sigma}(u)) = 0 & {\rm in \ } \Omega, \\
    \boldsymbol{\sigma}(u)  \cdot n = g^N & {\rm on \ } \Gamma^N, \\
    u = g^D & {\rm on \ } \Gamma^D,
  \end{cases}
\label{eq:elasticity_problem}
\end{equation}
with boundary data $g^D$ and $g^N$ set in accordance with the analytical solution of the infinite plate problem \eqref{eq:an_cutplate}. The Cauchy stress is related to the displacement field by Hooke's Law, $\boldsymbol{\sigma}(u) = \lambda  {\rm div}(u) \mathbf{I} + 2 \mu \nabla^s u$, with $\nabla^s$ the symmetric gradient operator. 

In view of symmetry, only a quarter of the domain with length $L=1$ is considered (Figure~\ref{fig:cutplate_solution}). The domain is rotated over an angle of \SI{45}{\degree} with respect to the Cartesian mesh with size $h$ in which it is immersed (Figure~\ref{fig:cutplate_grid}). Quadratic B-spline spaces are constructed over a sequence of ten uniformly refined meshes, starting with cells of size $h=1$. The number of elements that intersect the domain ranges from $4$ on the coarsest mesh to $216{,}672$ on the finest mesh. To preserve the geometry parametrization provided by the bisection-based tessellation scheme, the maximal recursion depth is reduced with each mesh refinement. On the finest mesh (nine refinements) trimmed cells are directly triangulated for integration purposes. Note that exact geometric representation of the domain would yield homogeneous Neumann data: $g^N = 0$. However, since the employed tessellation is only an approximation of the exact geometry, a (small) traction is applied along the circular boundary in accordance with \eqref{eq:an_cutplate} and the approximate geometry.

The finite dimensional weak formulation corresponding to \eqref{eq:elasticity_problem} is of the form \eqref{eq:weakformproblem} with (see \citep{Ruess2013} for details)
\begin{subequations}\label{eq:bilop_2d}\begin{align}
\mathcal{F}(v_h,u_h)  = & \int_\Omega \nabla^s v_h : \boldsymbol{\sigma}(u_h) {\rm d}V \nonumber \\
           & + \int_{\Gamma^D} - \big( v_h \cdot \boldsymbol{\sigma}(u_h) \cdot n + u_h \cdot \boldsymbol{\sigma}(v_h) \cdot n \big) {\rm d}S \nonumber \\
           & + \int_{\Gamma^D} \big( \beta^\lambda v_h \cdot (n \otimes n) \cdot u_h + \beta^\mu v_h \cdot u_h \big) {\rm d}S, \\
\ell(v)  = & \int_{\Gamma^N} v_h \cdot g^N {\rm d}S + \int_{\Gamma^D} -g^D \cdot \boldsymbol{\sigma}(v_h) \cdot n {\rm d}S \nonumber \\
           &  + \int_{\Gamma^D} \big( \beta^\lambda v_h \cdot (n \otimes n) \cdot g^D + \beta^\mu v_h \cdot g^D \big) {\rm d}S.
\end{align}\end{subequations}
Local stabilization with parameters $\beta^\lambda = 2 \lambda C^\lambda > \lambda C^\lambda$ and $\beta^\mu = 4 \mu C^\mu > 2 \mu C^\mu$ is used, with $C^\lambda$ and $C^\mu$ satisfying
\begin{subequations}\begin{align}
\int_{\Gamma^D} {\rm div} (v_h)^2 \text{d}\Gamma & \leq \int_\Omega C^\lambda {\rm div} (v_h)^2 \text{d}\Omega \quad \forall v_h \in \mathcal{V}_h, \\
\int_{\Gamma^D} (n \cdot \nabla^s v_h)^2 \text{d}\Gamma & \leq \int_\Omega C^\mu |\nabla^s v_h|^2 \text{d}\Omega \quad \forall v_h \in \mathcal{V}_h.
\end{align}\end{subequations}

\begin{figure}
  \centering
  \begin{subfigure}[t]{0.49\textwidth}
    \includegraphics[width=\textwidth]{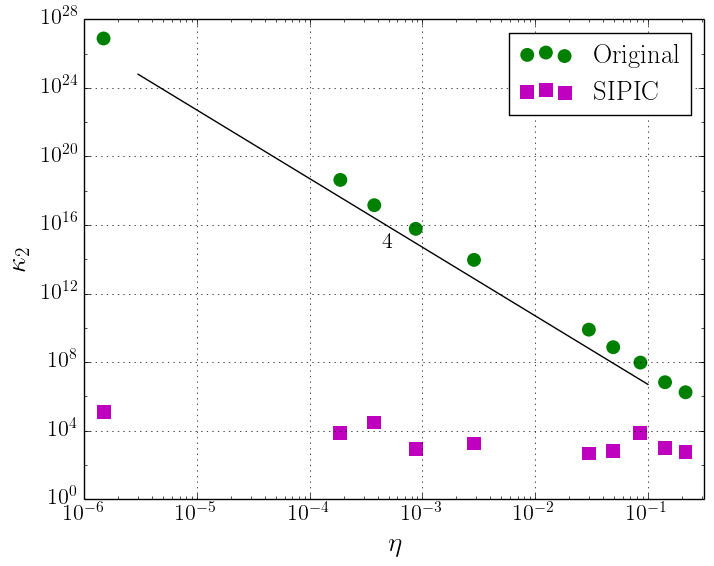}
    \caption{}
    \label{fig:cutplate_kappa_eta}
  \end{subfigure}
  \begin{subfigure}[t]{0.49\textwidth}
    \includegraphics[width=\textwidth]{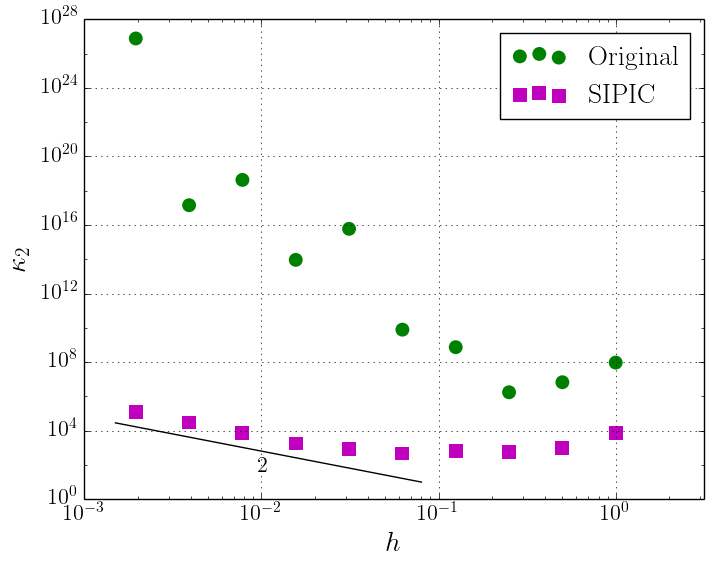}
    \caption{}
    \label{fig:cutplate_kappa_h}
  \end{subfigure}
  \caption{Original and SIPIC-preconditioned condition number \emph{vs.}\ {(a)} the smallest volume fraction and {(b)} the mesh size, for the plate with a circular hole with second order B-spline bases.}
  \label{fig:cutplate_kappa}
\end{figure}

Figure~\ref{fig:cutplate_kappa} shows the original and SIPIC-preconditioned condition numbers versus the smallest volume fraction (Figure~\ref{fig:cutplate_kappa_eta})
and the mesh size (Figure~\ref{fig:cutplate_kappa_h}). The scaling relation \eqref{eq:estimate} for the original system is clearly observed from Figure~\ref{fig:cutplate_kappa_eta}, despite the fact that the mesh size $h$ is varied. One can observe that the SIPIC-preconditioned condition number is essentially independent of the minimal volume fraction. Both observations are in agreement with the numerical results in Section~\ref{sec:preconexample}. Figure~\ref{fig:cutplate_kappa_h} shows the same data versus the mesh size $h$. The observed relation between the original condition number and the mesh size  can be conceived of as an induced effect, due to a correlation between the mesh size $h$ and smallest volume fraction $\eta$, \emph{viz.,} in general smaller minimum volume fractions can be expected when the number of trimmed cells is increased (by $h$-refinement). An important observation from Figure~\ref{fig:cutplate_kappa_h} is that SIPIC preconditioning results in a scaling relation $\kappa_2(\mathbf{SAS}^T) \propto h^{-2}$, which resembles that of standard finite elements \cite{johnson}.

\begin{figure}
  \centering
  \begin{subfigure}{0.49\textwidth}
    \includegraphics[width=\textwidth]{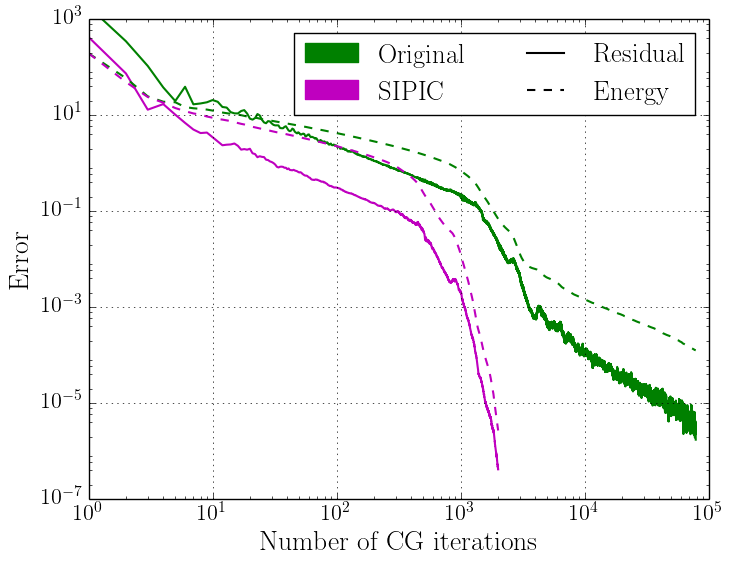}
    \caption{CG iterations}\label{fig:cg_convergence_cutplate}
  \end{subfigure}
  \begin{subfigure}{0.49\textwidth}
    \includegraphics[width=\textwidth]{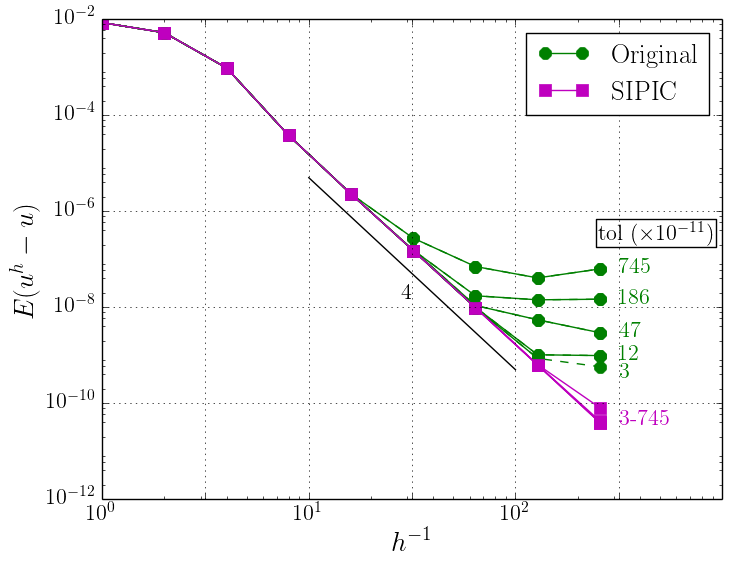}
    \caption{Mesh convergence}\label{fig:h_refinement_cutplate}
  \end{subfigure}
    \caption{Original and SIPIC-preconditioned finite cell results for the plate with a circular hole with second order B-spline bases: {(a)} Convergence of the residual and energy norm error \emph{vs.}\ the CG iterations on the finest mesh; {(b)} Convergence of the strain energy error under mesh refinement for various CG solver tolerances.}
    \label{fig:plateperformance}
\end{figure}

The performance of the SIPIC preconditioner is studied further in Figure~\ref{fig:plateperformance}. Figure~\ref{fig:cg_convergence_cutplate} displays the convergence of the CG solver for the finest mesh ($438{,}756$ degrees of freedom) with and without SIPIC preconditioning. The original condition number $\kappa_2(\mathbf{A}) \approx 10^{27}$, whereas the SIPIC-preconditioned condition number equals $\kappa_2(\mathbf{S}\mathbf{A}\mathbf{S}^Ts) \approx 10^5$. Both the residual ($\|\mathbf{b} - \mathbf{A}\mathbf{x}_i\|_2$, solid lines) and the error in the energy norm ($\|\mathbf{x} - \mathbf{x}_i\|_{\mathbf{A}}$, dashed lines) are shown. As expected from Equation~\eqref{eq:cgconvergence}, it is indeed observed that the SIPIC-preconditioned system converges significantly faster than the original system. For a CG solver tolerance of $10^{-6}$, the preconditioned system requires approximately $50$ times fewer iterations. Furthermore, following Equation~\eqref{eq:cgerrorbounds}, it is observed that the correlation between the residual and the energy norm error is stronger in the preconditioned system. For the same tolerance, the energy norm error in the original system is substantially larger than the energy norm error in the preconditioned system. 

An important effect of the stricter correlation between the residual and the energy norm error for preconditioned systems is visible in Figure~\ref{fig:h_refinement_cutplate}. In this figure, we study the mesh convergence by plotting the strain energy of the error, $E(u-u_h) = \frac{1}{2}\nabla^s(u-u_h):\boldsymbol{\sigma}(u-u_h)$, against the mesh size for various CG solver tolerances (ranging from $3\times 10^{-11}$ to $745\times 10^{-11}$). The CG solver terminates when either the relative residual ($\|\mathbf{b}-\mathbf{A}\mathbf{x}_i\|_2/\|\mathbf{b}\|_2$) or the absolute residual ($\|\mathbf{b}-\mathbf{A}\mathbf{x}_i\|_2$) reaches the specified tolerance, or when the number of CG iterations exceeds $100{,}000$. One can observe that for the original system, the error due to not solving the linear system with sufficient accuracy becomes dominant when the mesh is refined. For these ill-conditioned systems, the precision with which the linear system can be solved hinders mesh convergence of the strain energy error under mesh refinement. Improving the accuracy of the original system's solution by lowering the CG tolerance is not a practical solution to this problem, since the number of iterations required to reach convergence of the CG solver (and thereby the computational effort) increases dramatically. For the the finest mesh with the smallest considered tolerance, a converged result was not obtained within $100{,}000$ CG iterations (indicated by the dashed line in Figure~\ref{fig:h_refinement_cutplate}). The SIPIC preconditioner improves the quality of the solution measured in the strain energy of the error, and asymptotic convergence under mesh refinement is observed even for relatively large CG solver tolerances. The observed rate of mesh convergence of the strain energy error of $h^4$ resembles the optimal rate for standard finite elements of order $p=2$.

\subsection{Quarter solid torus}
\label{sec:quarter_ring}
\begin{figure}
 \centering
 \begin{subfigure}[t]{.53\textwidth}
  \includegraphics[width=\textwidth]{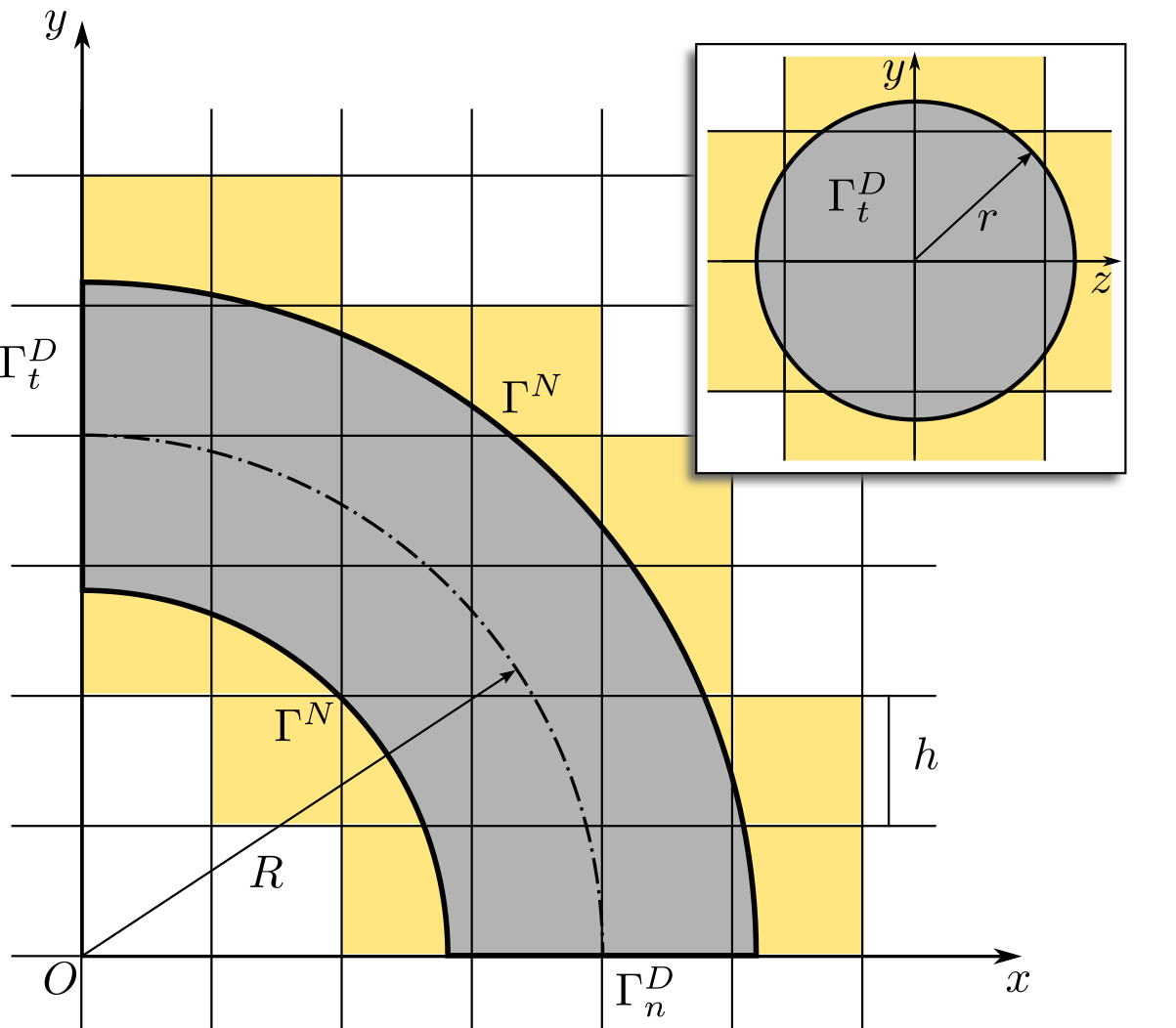}
  \caption{}
  \label{fig:quarter_ring_setup}
 \end{subfigure}
 \begin{subfigure}[t]{.46\textwidth}
  \includegraphics[width=\textwidth]{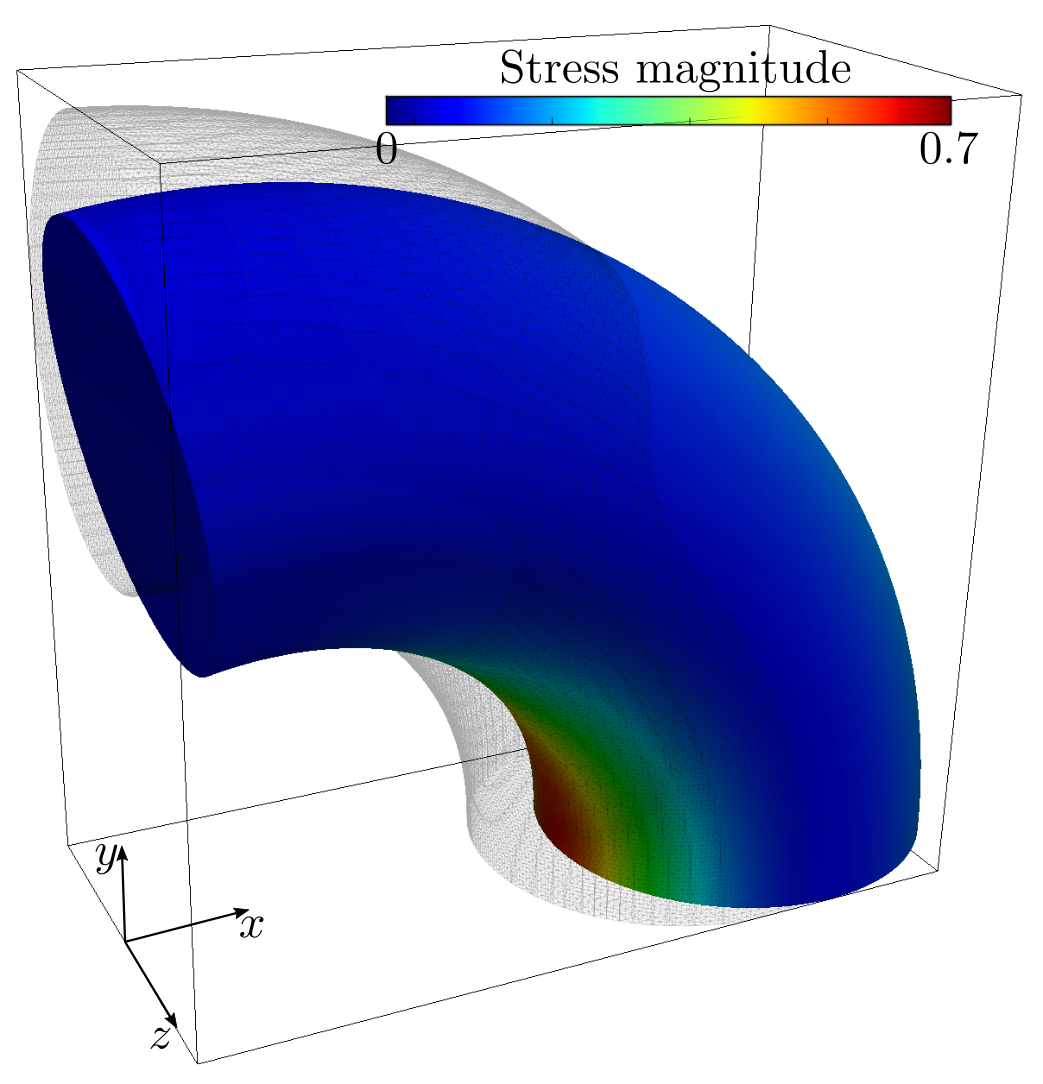}
  \caption{}
  \label{fig:quarter_ring_solution}
 \end{subfigure}
 \caption{Problem setup for a quarter of a solid torus with radius of revolution $R$ and cross-sectional radius $r$. The flat ends of the quarter torus coincide with the mesh. The stress magnitude is defined as the Frobenius norm of the Cauchy stress tensor, and the displacements in (b) are scaled.}\label{fig:quarter_ring_testcase}
\end{figure}

We consider the linear elastic analysis of a quarter of a three-dimensional solid torus (Figure~\ref{fig:quarter_ring_testcase}). The radius of the torus, measured from its axis of revolution ($z$-axis) to the center of its circular cross-sectional area, is equal to $R=2$. The radius of the cross-section is equal to $r=\frac{1}{2}\sqrt{2}$. On the bottom boundary ($\Gamma_n^D$, $y=0$) only the normal ($y$-)component of the displacement is constrained. On the left boundary ($\Gamma_t^D$, $x=0$) the downward ($y$-)displacement is prescribed, while the displacement in the $z$-direction is constrained. This boundary can move freely in the normal ($x$-)direction. 

The strong formulation is given by
\begin{equation}
     \begin{cases}
       -{\rm div} ( \boldsymbol{\sigma} ) = 0 & \mbox{in } \Omega, \\
       u \cdot n = 0  \quad \mbox{and} \quad \boldsymbol{P} \cdot \boldsymbol{\sigma} \cdot n = 0 & \mbox{on } \Gamma_n^D, \\
       \boldsymbol{P} \cdot u = \boldsymbol{P} \cdot g^D   \quad \mbox{and} \quad n \cdot \boldsymbol{\sigma} \cdot n = 0 & \mbox{on } \Gamma_t^D, \\
       \boldsymbol{\sigma} \cdot n = 0  & \mbox{on } \partial \Omega \setminus ( \Gamma_n^D \cup \Gamma_t^D ),
     \end{cases}
     \label{eq:elasticity_problem_3d}
\end{equation}
with projection tensor $\boldsymbol{P} = \boldsymbol{I} - n \otimes n$ and Dirichlet data $g^D = (0,0,-1)$. Note that this condition is only applied on the tangential (left) boundary $\Gamma_t^D$, as on the normal (bottom) boundary $\Gamma_n^D$ homogeneous Dirichlet conditions are applied. The Cauchy stress tensor is related to the displacement field by Hooke's law with Lam\'e parameters $\lambda = \mu = 1$. Note that this problem is underconstrained as it permits for a rigid body translation in the $x$-direction and an infinitesimal rigid body rotation around the $y$-axis. Herein, these rigid body modes are automatically accounted for by the employed CG solver.

The solid torus is immersed in a Cartesian mesh with size $h$ that is aligned with the coordinate system in Figure~\ref{fig:quarter_ring_setup}. Quadratic B-spline spaces are considered over a sequence of six uniformly refined meshes, starting with a mesh of size $h=1$. The recursion depth for the integration scheme is equal to five on the coarsest mesh, and is decreased upon mesh refinement to preserve the geometry parametrization. The number of elements that intersect the computational domain ranges from $16$ on the coarsest mesh to $168{,}928$ on the finest mesh.

The operators for the finite dimensional weak formulation \eqref{eq:weakformproblem} corresponding to problem \eqref{eq:elasticity_problem_3d} -- with the Dirichlet boundary conditions imposed by Nitsche's method -- are given by
\begin{subequations}\label{eq:bilop_3d}\begin{align}
\mathcal{F}(v_h,u_h) & = \int_\Omega \nabla^s v_h : \boldsymbol{\sigma}(u_h) {\rm d}V \nonumber \\
                     & + \int_{\Gamma_n^D} -\big( v_h \cdot (n \otimes n) \cdot \boldsymbol{\sigma}(u_h) \cdot n + u_h \cdot (n \otimes n) \cdot  \boldsymbol{\sigma}(v_h) \cdot n \big) {\rm d}S \nonumber \\
                     & + \int_{\Gamma_n^D} (\beta^\lambda + \beta_n^\mu) v_h \cdot (n \otimes n) \cdot u_h {\rm d}S \nonumber \\
                     & + \int_{\Gamma_t^D} -\big( v_h \cdot \boldsymbol{P} \cdot \boldsymbol{\sigma}(u_h) \cdot n + u_h \cdot \boldsymbol{P} \cdot \boldsymbol{\sigma}(v_h) \cdot n \big) {\rm d}S \nonumber \\
                     & + \int_{\Gamma_t^D} \beta_t^\mu v_h \cdot \boldsymbol{P} \cdot u_h {\rm d}S, \\
\ell(v)              & = \int_{\Gamma_t^D} - g^D \cdot \boldsymbol{P} \cdot \boldsymbol{\sigma}(v_h) \cdot n {\rm d}S \nonumber \\
                     & + \int_{\Gamma_t^D} \beta_t^\mu v_h \cdot \boldsymbol{P} \cdot g^D {\rm d}S.
\end{align}\end{subequations}
The local stabilization parameters are taken as $\beta^\lambda  = 2 \lambda C^\lambda > \lambda C^\lambda$, $\beta_n^\mu = 4 \mu C_n^\mu > 2 \mu C_n^\mu$, and $\beta_t^\mu = 4 \mu C_t^\mu > 2 \mu C_t^\mu$. Provided that each trimmed cell intersects with at most one type of Dirichlet boundary, the constants $C^\lambda$, $C_n^\mu$, and $C_t^\mu$ satisfy:
\begin{subequations}\begin{align}
\int_{\Gamma_n^D} {\rm div} (v_h)^2 \text{d}\Gamma & \leq \int_\Omega C^\lambda {\rm div} (v_h)^2 \text{d}\Omega \quad \forall v_h \in \mathcal{V}_h, \\
\int_{\Gamma_n^D} (n \cdot \nabla^s v_h \cdot n)^2 \text{d}\Gamma & \leq \int_\Omega C_n^\mu |\nabla^s v_h|^2 \text{d}\Omega \quad \forall v_h \in \mathcal{V}_h, \\
\int_{\Gamma_t^D} (\boldsymbol{P} \cdot \nabla^s v_h \cdot n)^2 \text{d}\Gamma & \leq \int_\Omega C_t^\mu |\nabla^s v_h|^2 \text{d}\Omega \quad \forall v_h \in \mathcal{V}_h.
\end{align}\label{eq:stab_3d}\end{subequations}
The weak imposition of Dirichlet constraints in normal or tangential direction by means of the above operators has, to the best of our knowledge, not been reported in literature. Since the Dirichlet boundaries are aligned with the Cartesian mesh, the Dirichlet constraints can also be imposed strongly by using open knot vectors in the parametric directions corresponding to the $x$- and $y$-axes. This enables us to compare weak imposition of boundary conditions using Nitsche's method to (traditional) strong imposition of boundary conditions.

\begin{figure}
  \centering
  \begin{subfigure}[t]{0.49\textwidth}
    \includegraphics[width=\textwidth]{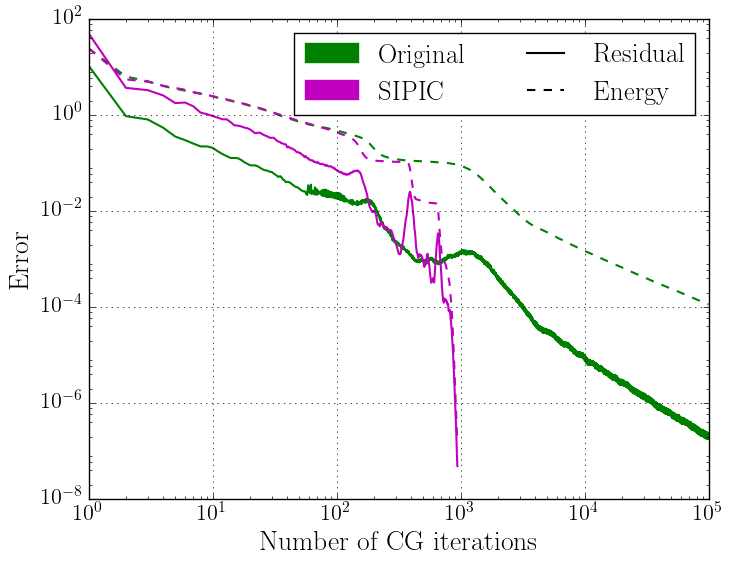}
    \caption{Nitsche BCs}
    \label{fig:cg_convergence_3d_weak}
  \end{subfigure}
  \begin{subfigure}[t]{0.49\textwidth}
    \includegraphics[width=\textwidth]{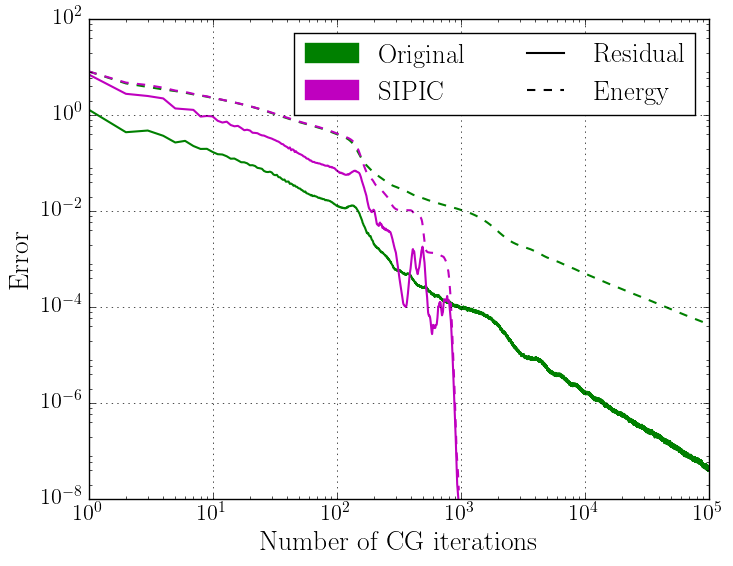}
    \caption{Strong BCs}
    \label{fig:cg_convergence_3d_strong}
  \end{subfigure}
  \caption{Convergence of the residual and energy norm error \emph{vs.}\ the solver iterations for the original and SIPIC-preconditioned finite cell results of the quarter solid torus problem with second order B-spline bases and a mesh size of $h=0.03125$ ($592{,}680$ degrees of freedom). Boundary conditions are {(a)} imposed using Nitsche's method, and {(b)} strongly imposed through the approximation space.}\label{fig:cg_convergence_3d}
\end{figure}

\begin{figure}
  \centering
  \includegraphics[width=.5\textwidth]{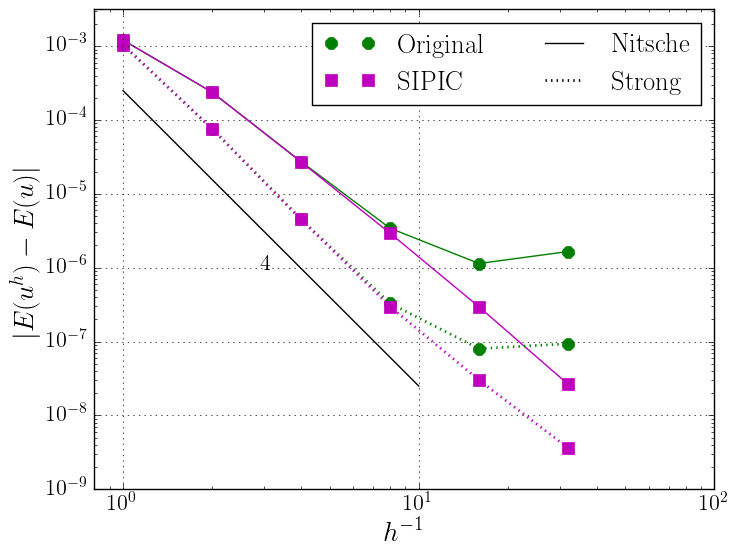}
  \caption{Convergence of the strain energy under mesh refinement for the original and SIPIC-preconditioned finite cell results for the quarter solid torus problem with second order B-spline bases. A CG solver tolerance of  of $9.5\cdot 10^{-7}$ is used.}\label{fig:h_refinement_3d}
\end{figure}

On the finest mesh, the number of degrees of freedom is equal to $592{,}680$, which makes the computation of condition numbers impractical. Therefore, we will here focus on the performance of the preconditioned iterative solver and the effect of SIPIC preconditioning on the mesh convergence behavior. Figure~\ref{fig:cg_convergence_3d} displays the CG solver convergence on the finest mesh, for both impositions of boundary conditions. One can observe that the convergence of the CG method is very similar for the formulations with weakly and strongly imposed boundary conditions, for both the original and the SIPIC-preconditioned system. Figure~\ref{fig:cg_convergence_3d} conveys a significant difference between the convergence behavior of CG for the original and the preconditioned system. SIPIC preconditioning reduces the number of CG iterations to reach a residual of $10^{-7}$ by more than a factor of $100$ relative to the original system. Moreover, SIPIC preconditioning improves the correlation between the residual ($\|\mathbf{b} - \mathbf{A}\mathbf{x}_i\|_2$) and the energy norm error ($\|\mathbf{x} - \mathbf{x}_i\|_{\mathbf{A}}$). The effect of the stricter correlation between the residual and the energy norm error for preconditioned systems is visible in Figure~\ref{fig:h_refinement_3d}, which studies mesh convergence of the error in the strain energy using a solver tolerance of $9.5\cdot 10^{-7}$. The error in the strain energy is defined as the absolute difference in strain energy compared to an overkill solution computed on the finest mesh with cubic B-splines (\emph{i.e.,} $|E(u_{\rm overkill}) - E(u_h)|$ with $E(u) = \frac{1}{2}\nabla^s u:\boldsymbol{\sigma}(u)$). Figure~\ref{fig:h_refinement_3d} clearly illustrates that ill-conditioning of the original system impedes asymptotic mesh convergence and hinders the computation of high-accuracy solutions. Using SIPIC preconditioning, optimal asymptotic mesh convergence behavior is observed. The rates of mesh convergence for both weakly and strongly imposed Dirichlet boundary conditions closely resemble the optimal rate of $h^4$ for standard finite elements, despite the fact that the mesh convergence rate can be affected by \emph{e.g.,} the mesh dependent bilinear forms considered herein, or by a lack of accuracy in the overkill solution.

\section{Concluding remarks}
\label{sec:conclusion}
In this work, we have rigorously established a scaling relation for the condition number of second-order elliptic finite cell systems. This scaling relation reveals that the condition number is inversely proportional to the smallest trimmed cell volume fraction, a result which confirms observations reported in the literature. Since this inverse proportionality scales with the power of $2p$, with $p$ the order of the considered basis functions, ill-conditioning problems become particularly apparent when higher-order methods (\emph{e.g.,} IGA or p-FEM) are considered. The obtained scaling relation is valid for locally stabilized systems, independent of the type of boundary conditions. For globally stabilized systems, an even stronger dependence of the condition number on the volume fraction is observed.

We have developed the SIPIC (\emph{Symmetric Incomplete Permuted Inverse Cholesky}) preconditioning technique, which effectively improves the conditioning of finite cell systems. This improved conditioning has been observed directly in condition number computations and indirectly in the iterative solver performance. Compared to alternative techniques for improving the conditioning of finite cell systems, an advantage of the SIPIC preconditioning technique is that it does not manipulate the weak formulation or the employed approximation space. On one hand, the preconditioner applies diagonal scaling to the applied basis to avoid ill-conditioning due to small (in the finite cell norm) basis functions. On the other hand, ill-conditioning due to quasi linearly dependent basis functions is mitigated by the application of local Gram-Schmidt orthonormalization. The preconditioner is constructed at negligible computational cost and is solely based on information from the finite cell system matrix, which makes it algebraic. By virtue of this algebraic nature, it is straightforward to implement as it is non-intrusive in the rest of the numerical code. Applying the preconditioner does not yield significant fill-in to the finite cell system matrix. For locally stabilized systems, we found that the condition number of the SIPIC-preconditioned system matrix is independent of the smallest trimmed cell volume fraction, which is essential to enable robust solutions of finite cell systems. Although it is observed that for globally stabilized formulations a slight dependence of the condition number on the smallest volume fraction remains, the preconditioning technique also dramatically improves the conditioning of such systems.

The SIPIC preconditioner is closely related to well-established preconditioning techniques, in particular diagonal scaling and incomplete Cholesky preconditioning. We anticipate that alternative Cholesky-based preconditioning strategies can also improve the conditioning of finite cell systems. An advantage of the SIPIC preconditioner is that it is completely targeted to ill-conditioning sources related to small volume fractions. A detailed comparative analysis with alternative, Cholesky-based, preconditioning techniques constitutes an intricate task and is beyond the scope of our current work. In this manuscript, the preconditioner is applied symmetrically, in the perception of preconditioning the basis of the approximation space rather than preconditioning the system matrix. However, the preconditioner can also be applied in standard preconditioned conjugate gradient algorithms.

Herein we have restricted ourselves to finite cell systems. However, the sources of ill-conditioning indicate that similar condition number scaling relations hold for other immersed methods, enriched finite element methods such as XFEM, and weak coupling strategies. We expect the SIPIC preconditioning technique to also be effective for such problems. In this study, we considered piecewise polynomial bases, but the SIPIC preconditioning technique may also be effective with other types of bases. A particular point of attention for the application of SIPIC to different methods or approximation spaces is the selection of an appropriate orthonormalization threshold.

The current work has been restricted to symmetric positive definite systems, such that the bilinear form is an inner product with induced norm. Evidently, this covers a large variety of problems in computational mechanics. However, problems such as convection-dominated flows and mixed problems such as (Navier-)Stokes flows and incompressible elasticity are not within the range of applicability of this work. Generalization of the developments in this manuscript to such problems is a topic of further study.

\section*{Acknowledgement}
\noindent

The research of F.\ de Prenter was funded by the NWO under the Graduate Programme Fluid \& Solid Mechanics. The research of C.V.\ Verhoosel was funded by the NWO under the VENI scheme. All simulations in this work were performed using the open source software package Nutils (www.nutils.org).

\bibliographystyle{plain}
\bibliography{reference}

\appendix
\renewcommand*{\thesection}{\appendixname~\Alph{section}}
\renewcommand*{\theequation}{\Alph{section}.\arabic{equation}}

\section{Local eigenvalue problem conditioning}
\label{sec:local}
Determination of the local stability parameter $\beta_i > C_i = \lambda_{\rm max}$, requires the computation of the largest eigenvalue of the local eigenvalue problem (see Section \ref{sec:coercivity}):
\begin{equation}
\mathbf{B}_i\mathbf{y} = \lambda \mathbf{V}_i\mathbf{y}.
\tag{\ref{eq:geneig}}
\end{equation}
The numerical solution of this problem requires careful consideration of the following two aspects:
\begin{itemize}
 \item The matrix $\mathbf{V}_i$, equation \eqref{eq:localevV}, is singular since it is based solely on the gradients of the shape functions in $\boldsymbol{\Phi}_i$. This singularity is not a fundamental problem for the computation of $C_i$ in accordance with \eqref{eq:infdim}, by virtue of the fact that ${\rm ker}(\mathbf{V}_i) \subset {\rm ker}(\mathbf{B}_i)$. However, the kernel of $\mathbf{V}_i$ can negatively affect the robustness with which the largest eigenvalue can be computed. 
 \item After removal of the kernel of $\mathbf{V}_i$, its conditioning can remain poor on cells with small volume fractions due to quasi linear dependence of basis functions (see Section~\ref{sec:orthonormalization}).
\end{itemize}
The conditioning of $\mathbf{V}_i$ of can be improved using the SIPIC preconditioner, but for these matrices a more direct solution exists. Since the eigenvalue problem is local to the cells, the numerical difficulties related to these conditioning aspects can be circumvented by applying a local change of basis. Quasi linear dependencies are avoided by employing the basis of monomials
\begin{equation}
 \phi_{\boldsymbol{\rho}} = \prod_{j=1}^d \left( x_j - \hat{x}_j \right)^{{\rho_j}}, 
\end{equation}
with multi-index ${\boldsymbol{\rho}}=(\rho_1,\ldots,\rho_d)$ for $0 \leq \rho_i \leq p$ and $\hat{x}$ the center of mass of the trimmed cell. With this basis, diagonal scaling suffices for the robust computation of $C_i$. Centering the monomials at the center of mass of the trimmed cell minimizes the quasi linear dependence problem because higher order terms occur in separate functions, such that after diagonal scaling their contribution is not diminished on small volume fractions of a cell. Furthermore, the kernel of $\mathbf{V}_i$ only consists of the function $\phi_{\boldsymbol{0}}$, such that removing the kernel is straightforward.

\section{Upper bounds for $\mathcal{F}^2(\cdot,\cdot)$ and $\mathcal{F}^3(\cdot,\cdot)$}\label{sec:app_bound}
Under the assumptions in Section~\ref{sec:conditioning}, we can bound $\mathcal{F}^2(\cdot,\cdot)$ and $\mathcal{F}^3(\cdot,\cdot)$ from above.
For the polynomial basis, we additionally assume that $\exists C_{\mathcal{V}}$ ($\sim 1$), independent of mesh size $h$, such that for every basis function $\phi$ in basis $\boldsymbol{\Phi}$ it holds that
\begin{equation}
 \|\phi\|_{L^\infty} \leq C_{\mathcal{V}} \qquad\mbox{and}\qquad \|\nabla \phi\|_{L^\infty} \leq \frac{C_{\mathcal{V}}}{h}.
\end{equation}
On cell $\Omega_i$ we can then bound $\mathcal{F}_i^2(\cdot,\cdot)$ and $\mathcal{F}_i^3(\cdot,\cdot)$:
\begin{subequations}\begin{align}
 \mathcal{F}_i^2(\phi,\phi) & \leq \int_{\Gamma_i^D} 2 \phi |\partial_n \phi| {\rm d}S  \leq 2 |\Gamma_i^D| \|\phi\|_{L^\infty} \| \nabla \phi \|_{L^\infty} \\
 & \leq 2 C_\Gamma C_{\mathcal{V}}^2 h^{d-2} \eta_i^{\frac{d-1}{d}}, \nonumber \\
 \mathcal{F}_i^3(\phi,\phi) & = \int_{\Gamma_i^D} \beta_i \phi^2 {\rm d}S  \leq |\Gamma_i^D| \beta_i \|\phi\|_{L^\infty}^2  \leq C_\Gamma C_\beta C_{\mathcal{V}}^2 h^{d-2} \eta_i^{\frac{d-2}{d}}. \label{eq:appcbound}
\end{align}\label{eq:bc_separate}\end{subequations}
For $d\geq2$, the inequalities in \eqref{eq:bc_separate} can be aggregated as
\begin{equation}
 \mathcal{F}_i^2(\phi,\phi) + \mathcal{F}_i^3(\phi,\phi) \leq C_{\mathcal{F}_i^{23}} h^{d-2} \eta_i^{\frac{d-2}{d}} \leq C_{\mathcal{F}_i^{23}} h^{d-2},
\end{equation}
with $C_{\mathcal{F}_i^{23}} = \left(2+C_\beta\right) C_\Gamma C_{\mathcal{V}}^2$. A similar derivation conveys that for the cross term between two different basis functions $\phi_\alpha$ and $\phi_\beta$
\begin{equation}
 \mathcal{F}_i^2(\phi_\alpha,\phi_\beta) + \mathcal{F}_i^3(\phi_\alpha,\phi_\beta) \leq C_{\mathcal{F}_i^{23}} h^{d-2}.
\end{equation}
The number of functions that is supported on a cell is $(p+1)^{d}$ for all piecewise polynomial spaces and the number of cells that support a function is is at most $(p+1)^d$ for B-splines and $2^d$ for most other bases.
This bounds the sum of all terms per function by $C_{\mathcal{F}^{23}}h^{d-2}$, with $C_{\mathcal{F}^{23}}^{\rm spline} = (p+1)^{2d}C_{\mathcal{F}_i^{23}}$ for B-splines and by $C_{\mathcal{F}^{23}}^{\rm other} =  (p+1)^{d}2^dC_{\mathcal{F}_i^{23}}$ for most other bases. As a result:
\begin{equation}
 \|\begin{bmatrix}\mathcal{F}^2(\boldsymbol{\Phi},\boldsymbol{\Phi}^T) + \mathcal{F}^3(\boldsymbol{\Phi},\boldsymbol{\Phi}^T)\end{bmatrix}\mathbf{y}\|_2 \leq C_{\mathcal{F}^{23}}h^{d-2} \|\mathbf{y}\|_2,
\end{equation}
such that
\begin{equation}
\max_{v_h,\mathbf{y}} \frac{\mathcal{F}^2(v_h,v_h)+\mathcal{F}^3(v_h,v_h)}{\mathbf{y}^T\mathbf{y}} \leq C_{\mathcal{F}^{23}}h^{d-2}.
\label{eq:noholder}
\end{equation}
The bounds in \eqref{eq:bcbounds} follow directly from \eqref{eq:noholder}. Inequality \eqref{eq:noholder} also follows directly from H\"{o}lders inequality. Because the sum of all terms per function is bounded by $C_{\mathcal{F}^{23}}h^{d-2}$, this value also bounds the $\|\cdot\|_\infty$ norm and the $\|\cdot\|_1$ norm (due to symmetry) of matrix $\begin{bmatrix}\mathcal{F}^2(\boldsymbol{\Phi},\boldsymbol{\Phi}^T) + \mathcal{F}^3(\boldsymbol{\Phi},\boldsymbol{\Phi}^T)\end{bmatrix}$. By H\"{o}lders inequality this implies $\|\begin{bmatrix}\mathcal{F}^2(\boldsymbol{\Phi},\boldsymbol{\Phi}^T) + \mathcal{F}^3(\boldsymbol{\Phi},\boldsymbol{\Phi}^T)\end{bmatrix}\|_2 \leq C_{\mathcal{F}^{23}}h^{d-2}$.

When a global stabilization parameter is applied, the bound for $\beta_i$ does not depend on $\eta_i$ but on $\eta$. The bound in \eqref{eq:appcbound} is then replaced by 
\begin{equation}
  \mathcal{F}_i^3(\phi,\phi) \leq |\Gamma_i^D| \beta \|\phi\|_{L^\infty}^2 \leq C_\Gamma C_\beta C_{\mathcal{V}}^2 h^{d-2} \eta_i^{\frac{d-1}{d}}\eta^{\frac{-1}{d}}.
\end{equation}
Assuming that the bound is attained, the following inequality holds:
\begin{equation}
 \max_{v_h,\mathbf{y}} \frac{\mathcal{F}^3(v^h,v^h)}{\mathbf{y}^T\mathbf{y}} \geq C h^{d-2} \eta^{\frac{-1}{d}}.
\end{equation}
Therefore, the condition number of globally stabilized systems is estimated by
\begin{equation}
\kappa_2(\mathbf{A}) \geq C \eta^{-(2p+1-1/d)},
\label{eq:estimate_global}
\end{equation}
which is different from \eqref{eq:estimate} and explains the different slope that is observed when global stabilization is applied in Section~\ref{sec:cond_example} and \ref{sec:preconexample}.

\section{Machine precision instabilities in the SIPIC construction algorithm}
\label{sec:app_singular}
The machine precision stability problem in the SIPIC construction algorithm as discussed in Section~\ref{sec:algorithm} is here illustrated for the $2\times 2$ matrix
\begin{align}
 \mathbf{A} = \begin{bmatrix} 1 & 1 - \varepsilon^2 \\ 1 - \varepsilon^2 & 1 \end{bmatrix},
 \label{eq:smallAmatrix}
\end{align}
with some parameter $\varepsilon \ll 1$ and condition number $\kappa_2 = 2/\varepsilon^2 - 1 \gg 1$. The SIPIC preconditioner for $\mathbf{A}$ can be computed analytically as
\begin{equation}
 \mathbf{S} = \begin{bmatrix} 1 & 0 \\ \frac{\varepsilon^2 - 1}{\varepsilon \sqrt{2-\varepsilon^2}} & \frac{1}{\varepsilon \sqrt{2-\varepsilon^2}} \end{bmatrix}.
\end{equation}
The computation of this matrix by the SIPIC construction algorithm can be impeded by machine precision errors when $\varepsilon^2 \sim {\tt eps}$. From \eqref{eq:smallAmatrix}, it is directly observed that in this case the rows of matrix $\mathbf{A}$ become linearly dependent up to machine precision. In the preconditioner construction algorithm, this results in the division by zero in the renormalization step after application of the Gram-Schmidt procedure. For this particular case, the Gram-Schmidt matrix equals 
\begin{equation}
 \mathbf{G} = \begin{bmatrix} 1 & 0 \\ \varepsilon^2 - 1 & 1 \end{bmatrix},
\end{equation}
which results in the orthogonalized matrix $\mathbf{A}^\perp=\mathbf{G}\mathbf{A} \mathbf{G}^T$ (since $\mathbf{D}={\rm diag}(\mathbf{A})=\mathbf{I}$). Normalization of this matrix results in the division by the square root of
\begin{equation}\begin{aligned}
A^\perp_{22}  &= \begin{bmatrix} \varepsilon^2 - 1 & 1 \end{bmatrix} \begin{bmatrix} 1 & 1 - \varepsilon^2 \\ 1 - \varepsilon^2 & 1 \end{bmatrix} \begin{bmatrix} \varepsilon^2 - 1 \\ 1 \end{bmatrix} \\
 & = 1^2 - (1-\varepsilon^2)^2,
\end{aligned}\end{equation}
which numerically becomes equal to zero when $\varepsilon^2 \sim {\tt eps}$.

\end{document}